\documentclass[preprint, 10pt]{elsarticle}

\usepackage{amssymb}
\usepackage{algorithm}
\usepackage{algorithmic}
\usepackage{tikz}
\usepackage{pgfplots}
\usepackage{microtype}
\usepackage{amsthm}
\usepackage{caption}

\usepackage{amsbsy}
\usepackage{amsmath}
\usepackage{bm}
\usepackage{color}
\usepackage{dcolumn}
\usepackage{tabularx}
\usepackage{bbm}

\newcommand{\vx}{{\mathbf x}}
\newcommand{\vy}{{\mathbf y}}
\newcommand{\vz}{{\mathbf z}}
\newcommand{\vv}{{\mathbf v}}
\newcommand{\vw}{{\mathbf w}}

\newcommand{\vp}{{\mathbf p}}
\newcommand{\vb}{{\mathbf b}}

\newif\ifnotesw \noteswtrue
\newcommand{\jeffrey}[1]{\ifnotesw  \textcolor[rgb]{0,0,1}{  $\spadesuit$Jeffrey:\ {\sf \bf \it #1}\ $\spadesuit$ }\fi}

\newcommand{\johannes}[1]{\ifnotesw  \textcolor[rgb]{1,0,1}{ $\square$ Johannes : {\sf \bf \it #1}\ $\square$ }\fi}
\noteswfalse   

\newcommand{\setalglineno}[1]{%
	\setcounter{ALC@line}{\numexpr#1-1}}
\newtheorem{theorem}{Theorem}[]

\newtheorem{lemma}[theorem]{Lemma}
\theoremstyle{definition}

\theoremstyle{remark}
\newtheorem*{remark}{Remark}

\journal{Journal of Computational and Applied Mathematics}

\begin{document}

\begin{frontmatter}



\title{Variable resolution Poisson-disk sampling for meshing discrete fracture networks}


\author{Johannes Krotz\corref{cor1} \fnref{label1}}

\affiliation[label1]{organization={Department of Mathematics, University of Tennessee},
            city={Knoxville},
            postcode={37919}, 
            state={Tennessee},
            country={USA}}

\author{Matthew R. Sweeney\fnref{label2}}

\author{Carl W. Gable\fnref{label2}}

\author{Jeffrey D. Hyman\fnref{label2}}

\affiliation[label2]{organization={Computational Earth Science (EES-16), Earth and Environmental Sciences, Los Alamos National Laboratory},
            city={Los Alamos},
            postcode={87545}, 
            state={New Mexico},
            country={USA}}

\author{Juan M. Restrepo\fnref{label3}}

\affiliation[label3]{organization={Oak Ridge National Laboratory},
            city={Oak Ridge},
            postcode={37830}, 
            state={Tennessee},
            country={USA}}

\begin{abstract}
We present the near-Maximal Algorithm for Poisson-disk Sampling (nMAPS) to generate point distributions for variable resolution Delaunay triangular and tetrahedral meshes in two and three-dimensions, respectively.
nMAPS consists of two principal stages.
In the first stage, an initial point distribution is produced using a cell-based rejection algorithm. 
In the second stage, holes in the sample are detected using an efficient background grid and filled in to obtain a near-maximal covering. 
Extensive testing shows that nMAPS generates a variable resolution mesh in linear run time with the number of accepted points.
We demonstrate nMAPS capabilities by meshing three-dimensional discrete fracture networks (DFN) and the surrounding volume.
The discretized boundaries of the fractures, which are represented as planar polygons, are used as the seed of 2D-nMAPS to produce a conforming Delaunay triangulation.
The combined mesh of the DFN is used as the seed for 3D-nMAPS, which produces conforming Delaunay tetrahedra surrounding the network.
Under a set of conditions that naturally arise in maximal Poisson-disk samples and are satisfied by nMAPS, the two-dimensional Delaunay triangulations are guaranteed to only have well-behaved triangular faces. 
While nMAPS does not provide triangulation quality bounds in more than two dimensions, we found that low-quality tetrahedra in 3D are infrequent, can be readily detected and removed, and a high quality balanced mesh is produced.

\end{abstract}



\begin{keyword}
Discrete Fracture Network\sep	Maximal Poisson-disk Sampling \sep Mesh Generation \sep Conforming Delauany Triangulation   



\end{keyword}

\end{frontmatter}


\section{Introduction}

There are a number of methods used to model flow and the associated transport of chemical species in low-permeability fractured rock, such as shale and granite.
The most common are continuum models, which use effective medium parameters~\cite{gerke1993dual,lichtner2014modeling,neuman1988use,neuman2005trends,tsang1996tracer,zimmerman1993numerical} and discrete fracture network/matrix (DFN) models, where fractures and the networks they form are explicitly represented~\cite{cacas1990modeling,long1982porous,nordqvist1992variable}. 
In the DFN methodology, individual fractures are represented as planar $N-1$ dimensional objects embedded within an $N$ dimensional space.
Each fracture in the network is resolved with a computational mesh and the governing equations for flow and transport are solved thereon.
While the explicit representation of fractures allows DFN models to represent a wider range of transport phenomena and makes them a preferred choice when linking network attributes to flow properties compared to continuum methods~\cite{hadgu2017comparative,hyman2018dispersion,hyman2019linking}, it also leads to unique and complex issues associated with mesh generation.
Both conforming methods, where the mesh conforms to intersections~\cite{hyman2014conforming,mustapha2007new,mustapha2011efficient}, and non-conforming methods, which use more complex discretization schemes so the mesh does not need to be conforming~\cite{berrone2013pde,erhel2009flow,pichot2010mixed,pichot2012generalized}, have been developed.  
If the volume representing the rock matrix surrounding the fracture network also needs to be meshed, then the complications associated with mesh generation are compounded for both techniques~\cite{berre2018flow}.

We present the near-Maximal Algorithm for Poisson-disk Sampling (nMAPS) that generates a conforming variable resolution triangulation mesh on and around three-dimensional discrete fracture networks.  
Although the theory of Poisson disk sampling to produce a point distribution is straightforward, details of the implementation, algorithm, termination criteria, and maximizing efficiency are nuanced and often problem-dependent. 
nMAPS addresses these issues in the context of mesh generation using a two-stage approach. 
The first stage is based on the framework presented by Dwork et al., \cite{fastvar} and uses a rejection algorithm to generate an initial Poisson-disk sample point set.
The second phase is based on the framework presented by Mitchell et al., \cite{varrad} and adds additional points to fill gaps in the covering, which maximizes point density without violating the restrictions of a Poisson-disk sampling.
nMAPS combines the two methods to achieve near-maximal coverage with a run time that scales linearly in the number of accepted points.
nMAPS efficiency results from a novel rejection technique that uses a quadrilateral or hexahedral background mesh search based on already sampled points.
This method significantly shortens the time needed to reach a sufficient, i.e., near-maximal, point density due to a large reduction in computational operations.
Once the point set is generated, the conforming Delaunay algorithm presented by Murphy et al.,~\cite{murphy2001point} is used to generate a Delaunay triangulation in two or three dimensions. 
Although nMAPS does not guarantee that the point distribution's density is maximal, i.e., no additional points can be added without violating the restrictions on distances between points, the sample sets are sufficiently maximal to produce high-quality meshes.
Comparison with previously implemented meshing techniques for variable resolution conforming triangulations showed that nMAPS produces meshes comparable in quality but requires considerably less computational time.  
Moreover, the nMAPS framework is presented in a general manner that can be easily extended beyond mesh generation for DFNs. 

In section~\ref{sec:background}, we describe the challenges associated with DFN mesh generation and the general properties of maximal Poisson-disk sampling. 
In section~\ref{sec:methods}, we provide a detailed explanation of nMAPS, for both 2D fracture networks and 3D volume meshing. 
In section~\ref{sec:results}, we assess the quality of the mesh and performance, i.e., run times, of 2D and 3D examples using nMAPS.
In section~\ref{sec:conclusions}, we provide a few concluding remarks.

\section{Background}\label{sec:background}

\subsection{Discrete Fracture Networks \& Mesh Generation}\label{ssec:dfn}

Due to the epistemic uncertainty associated with hydraulic and structural properties of subsurface fractured media, fracture network models are typically constructed stochastically~\cite{national2020characterization,national1996rock,neuman2005trends}.
In the DFN methodology, individual fractures are placed into the computational domain with locations, sizes, and orientations that are sampled from appropriate distributions based on field site characterizations.
The fractures form an interconnected network embedded within the porous medium. 
Each fracture must be meshed for computation so that the governing equations for flow and transport can be numerically integrated to simulate physical phenomena of interest.

Formally, each fracture in a DFN can be represented as a planar straight-line graph (PSLG) composed of a set of line segments that represent the boundary of the fracture and a set of line segments that represent where other fractures intersect it. 
In this manner, each fracture can be described by a set of boundary points on the PSLG, denoted $\{p\}$, and a set of intersection lines $\{\ell_{i,j}\}$, where the subscripts $i$ and $j$ indicate that this line corresponds to the intersection between the $i$th and $j$th fractures. 
Once $\{p\}$ and $\{\ell_{i,j}\}$ are obtained for every fracture in the network, a point distribution covering each fracture must be generated.
If a conforming numerical scheme is used, then all cells of $\{\ell_{i,j}\}$ are discretized lines in the mesh that must coincide between intersecting fractures. 
So long as minimum feature size constraints are met, a conforming triangulation method, such as that presented in Murphy et al.,~\cite{murphy2001point}, can be implemented to connect the vertices such that all lines of intersection form a set of connected edges in a triangulation.

In general, one wants to properly resolve all relevant flow and transport properties of interest while minimizing the number of nodes in the mesh, and these two goals compete.
The first of these conditions depends partially on the mesh quality, which is controlled partially by the second.
A starting point for the notion of mesh quality is the minimum angle condition, i.e., that the smallest angle should be bounded away from zero, as poorly shaped elements can affect the condition number of the linearized system of equations~\cite{azis,strangfix,ZLAMAL}. 
While a uniform resolution mesh is straightforward to generate and the minimum angle condition easily met, it is computationally more expensive than a variable resolution mesh refined in areas of interest, which can be tailored to reduce the number of nodes in the mesh.
A variable resolution mesh can be appropriate for single-phase flow simulations or in particle tracking simulations where the spatially variable resolution does not adversely affect transport properties.
Specifically, Eulerian formulations of transport will have spatially variable numerical diffusion on a variable resolution mesh, which means that the mesh needs to have a sufficiently small refinement to accurately capture solute fronts~\cite{benson2017comparison}.
Variable resolution mesh generation is more complex than its uniform counterpart. 
One of the principal complications is creating a smooth transition of cell sizes. 
Jumps in the computed fields of interest and other numerical artifacts can occur if the transition is not sufficiently smooth.

Most methods to generate a conforming DFN mesh use a uniform point distribution and modify connectivity locally to conform to intersections~\cite{mustapha2007new,mustapha2011efficient}.
When using a conforming mesh, the numerical methods for resolving flow and transport in the network are typically simpler and have fewer degrees of freedom compared to non-conforming mesh methods~\cite{fumagalli2019conforming}.  
Similarly, almost all non-conforming numerical methods use a uniform resolution, but some create variable resolutions across fractures, still uniform within a single plane, in an attempt to reduce the number of total nodes in the mesh~\cite{berrone2019parallel}.
Using a variable mesh resolution in non-conforming schemes could drastically reduce the number of nodes in the mesh while retaining the the ability to retain a high order of accuracy. 
However, it is rarely implemented due to the associated meshing complications~\cite{borio2021comparison}.

The generation of a variable resolution, unstructured conforming mesh is quite rare, even with the advantages noted above. 
One technique in use is the Features Rejection Algorithm for Meshing ({\sc FRAM}) that addressed the issues associated with conforming DFN mesh creation by coupling it with network generation~\cite{hyman2014conforming}. 
Through this technique, {\sc FRAM} allows for the creation of a variable resolution mesh that smoothly coarsens away from intersections where pressure gradients in flow simulations are typically the highest. 
{\sc FRAM} has been implemented in the computational suite {\sc dfnWorks}~\cite{hyman2015dfnWorks}, which has been used to probe fundamental aspects of geophysical flows and transport in fractured media~\cite{hyman2020flow,hyman2019emergence,hyman2019matrix,kang2020anomalous,makedonska2016evaluating,sherman2020characterizing} as well as practical applications including hydraulic fracturing operations~\cite{hyman2017discontinuities,karra2015effect,lovell2018extracting}, inversion of micro-seismicity data for characterization of fracture properties~\cite{mudunuru2017sequential}, the long term storage of spent civilian nuclear fuel~\cite{hadgu2017comparative}, and geo-sequestration of carbon dioxide into depleted reservoirs~\cite{hyman2020characterizing}.
However, the current method for mesh generation uses an inefficient iterative refinement technique to produce the point distribution used to generate the mesh. 
Sweeps of a refinement algorithm are applied to an initially coarse triangulation based on the boundary set $\{p\}$.
If an edge in the mesh is larger than the current maximum edge length, then a new point is added to the mesh at the midpoint of that edge to split it into two new edges. 
In practice, the edge splitting is done using Rivara refinement~\cite{rivara1984algorithms,rivara1984mesh}.
The resulting field is then smoothed using Laplacian smoothing in combination with Lawson flipping~\cite{khamayseh1996anisotropic}. 
This process is repeated until all edges meet the assigned target edge length, which could be a spatially variable field based on the distance to  $\{\ell_{i,j}\}$, for example.
While the resulting mesh quality is quite good, the process is rather slow and cumbersome.

If these implementation complexities can be addressed, then the superior modeling qualities of variable resolutions can be made practical.
To accomplish this, we present nMAPS, where the final vertex distribution is directly created rather than iteratively derived. 
While the method was initially designed to specifically improve {\sc FRAM}, we provide the details in a general format such that it can be implemented for mesh generation in general.
Specifically, it can be used as the basis for any DFN meshing methodology,  including both conforming and non-conforming techniques. 
Details are given for Delaunay triangulations, which are of importance in many two-point flux finite volume solvers as they are used to generate the Voronoi control volumes on which these solvers compute.  
In the next section, we recount the properties of maximal Poisson-disk sampling that we used to design and implement nMAPS, including theoretical bounds on mesh gradation that ensure high-quality variable mesh resolutions.

\subsection{Maximal Poisson-disk Sampling}\label{sec:mpds}

Meshes produced from a point distribution that is dense yet cluster free have provable high quality bounds under specific conditions~\cite{MDS_good_trig2,MDS_good_trig3,MDS_good_trig1}.
Maximal Poisson-disk samples fulfill these conditions~\cite{monte-carlo}. 
Similar quality bounds can be established for sphere-packings, whose radii are Lipschitz continuous with respect to their location~\cite{Lipschitz1,Lipschitz2,Lipschitz3}. 
Traditionally, Poisson-disk sampling is performed using an expensive dart-throwing algorithm~\cite{aliasing_3}.
These algorithms struggle to achieve maximality as the probability to select a free spot becomes decreasingly small as the number of points sampled ($n$) increases. 
The algorithm in \cite{varrad}, which is based on these dart-throwing algorithms, was the first to guarantee maximality. 
It achieves maximality with run times of $O(n\log(n))$ by using a regular grid for acceleration and sampling from polygonal regions in its second phase.
In practice, performance close to $O(n)$ has been reported~\cite{eff_mds,mds_simple,varrad}. 
Another algorithm that is not based on dart-throwing was proposed in \cite{Bridson2007FastPD}.
While the method does not guarantee maximality, it did show linear performance in the number of points sampled. 
This algorithm was extended to variable radii by \citep{fastvar}. 
Other authors have provided algorithms that produce variable Poisson-disk samplings on 3D-surfaces~\cite{GUO1,Guo2}. 
Their triangulation-based algorithm runs in linear time, and while theoretically not guaranteeing maximality, their numerical experiments indicate that near-maximality can certainly be achieved.  
A summary of recent developments in this area, along with a comparison of different methods, can be found in \cite{DongMing}.

Maximal Poisson-disk samplings $X$ on a domain $\Omega\subseteq \mathbb{R}^d$ are random selections of points $X=\{\vx_i\}_{i=1}^n$, that fulfill the following properties:
\begin{enumerate}
	\item \textit{empty disk property}: $$\forall i\neq j \in \{1,...,n\}: |\vx_i-\vx_j|>r.$$
	We will call $r$ the \textit{inhibition radius},
	\item \textit{maximality}:	
	$$\Omega= \bigcup_{i=1}^n B_R(\vx_i),$$
	where $B_\varepsilon(\vx)=\{\vy\in \Omega: |\vx-\vy|<\vy\}$ is the open ball of radius $\varepsilon$ around $\vx$. $R$ will be called the \textit{coverage radius}. \cite{varrad}
	
\end{enumerate}
Intuitively, the \textit{empty disk property} says that every sample point is at the center of a $d$-dimensional ball or disk that does not contain any other points of the sampling. 
\textit{Maximality} implies that these balls cover the whole domain, i.e., there is no point $y\in \Omega$ that is not already contained in one of the balls around a point in the sample.

It is useful to generalize these definitions, such that both the \textit{inhibition} and the \textit{coverage radius} depend on the sampling points, i.e., $r=r(\vx_i,\vx_j)$ and $R=R(\vx_i,\vx_j)$ for all $\vx_i,\vx_j\in X$. 
We hereafter refer to this construct as a variable radius maximal Poisson-disk sampling, and we refer to a Poisson-disk sampling with constant radii as a fixed-radii maximal Poisson-disk sampling~\cite{varrad}.

A common approach is to assign each point $\vx\in \Omega$ a positive radius $\rho(\vx)$ and have $r(\vx_i,\vx_j)$ be a function of $\rho(\vx_i)$ and $\rho(\vx_j)$. 
Natural choices for $r(\vx_i,\vx_j)$ are, for example, $\rho(\vx_i)$ or $\rho(\vx_j)$ for $i<j$, thereby determining the inhibition radius depending on the ordering on $X$. Order independent options include $\min(\rho(\vx_i),\rho(\vx_j)), \max(\rho(\vx_i),\rho(\vx_j))$ or $\rho(\vx_i)+\rho(\vx_j)$. 
The last of these options corresponds to a sphere packing \cite{varrad}. 
The coverage radius can but does not have to be different from $\rho$.


Delaunay triangulations maximize the smallest angle of all triangulations based on a generating point distribution \cite{def_delaunay}. 
Since numerical errors in many applications tend to increase if these angles are small \cite{ZLAMAL}, Delaunay triangulations often are a triangulation of choice.
Moreover, the dual of the Delaunay triangulation is a Voronoi tessellation, which in a certain sense is optimal for two-point flux finite volume solvers~\cite{eymard2000finite} that are commonly used in subsurface flow and transport simulators such as {\sc fehm}~\cite{zyvoloski2007fehm}, {\sc tough2}~\cite{pruess1999tough2}, and {\sc pflotran}~\cite{lichtner2015pflotran}. 
In the case of maximal Poisson-disk samplings, we can go one step further and give a lower bound on these angles. 
We provide a brief summary of the proofs found in \cite{varrad} and highlight the most important results we use. 
We then proceed with the new bounds.

\begin{lemma}
	\label{lemma:angle2d}
The smallest angle $\alpha$ in any triangle is greater than $\arcsin\left(\frac{r}{2R}\right)$, where $r$ is the length of the shortest edge and $R$ the radius of the circumcircle or \begin{align}
	\sin(\alpha)\ge \frac{r}{2R} \label{eq:smallest angle}
\end{align}
	\begin{proof}
		This is a direct corollary of the central angle theorem.
	\end{proof}
\end{lemma}
This Lemma allows us to give explicit bounds for maximal Poisson-disk samplings.
While we will focus entirely on inhibition radii given by $r(\vx_i,\vx_j)=\min(\rho(\vx_i),\rho(\vx_j))$, where $\rho(\vx)$ is some positive function, comparable results can be found for different $r(\vx_i,\vx_j)$ in a similar fashion. 
\begin{lemma}\label{lemma:2}
	Let $\varepsilon\ge 0$ and $\rho:\mathbb{R}^d\rightarrow \mathbb{R}$ ($d\ge 2$) be a positive Lipschitz continuous function with Lipschitz constant $L$ with $L\varepsilon < 1$.
	Let $X$ \jeffrey{I removed a $\subset$ please check that. There was no subset container. 
    Meaning that $X$ was not a subset of something else. Seemed like a typo.}\johannes{used to be $X\subset \Omega$, but felt redundant.}  be a variable maximal Poisson-disk sampling on the domain $\Omega\subset \mathbb{R}^d$ with inhibition radius $r(\vx,\vy)=\min(\rho(\vx),\rho(\vy))$  and coverage radius $R(\vx,\vy) \le (1+\varepsilon)r(\vx,\vy)$ with $\varepsilon>0$.
	Let the triangle $\Delta$ be an arbitrary element of the Delaunay triangulation of $X$ ($d=2$) or an arbitrary 2-dimensional face of a cell of the Delaunay triangulation of $X$ ($d > 2$).
    If the circumcenter of $\Delta$ is contained in $\Omega$, then each angle $\alpha$ of $\Delta$ is greater or equal to $\arcsin\left(\frac{1-L-\varepsilon L}{2+2\varepsilon}\right)$ or \jeffrey{should the second $\varepsilon$ be in the denominator, as in the following equation? I moved it there.  Check the previous see to see where it was. Please confirm}\johannes{yes}
	 \begin{align*}		 
	 \sin(\alpha)\ge \frac{1-L-\varepsilon L}{2+2\varepsilon}.
	\end{align*}

\begin{proof}
Let $\alpha$ be the smallest angle of $\Delta$ and $\vx,\vy\in X$ be the vertices of the shortest edge of $\Delta$, i.e., the vertices opposite to $\alpha$. 
Without loss of generality, assume $\rho(\vx)\le \rho(\vy)$.
Since $X$ is a Poisson-disk sampling $|\vx-\vy|\ge \min\left(\rho(\vx),\rho(\vy)\right)=\rho(\vx)$. 
Now, let $\vz\in \Omega$ be the circumcenter of $\Delta$. 
Since $X$ is maximal, there exists $\vv\in X$ with $|\vz-\vv|\le R(\vz,\vv)\le (1+\varepsilon)\rho(\vz)$.  
Because $\Delta$ was retrieved from a Delaunay triangulation $\vv$ cannot be contained in the interior of $\Delta$'s circumcircle. 
Hence,
\begin{align*}
	|\vz-\vx|\le |\vz-\vv|\le (1+\varepsilon)\rho(\vz)\le (1+\varepsilon)\left(\rho(\vx)+L|\vz-\vx| \right).
\end{align*}
Rearranging this inequality yields 
\begin{align*}
    |\vz-\vx|\le \rho(\vx)\frac{1+\varepsilon}{1-L-\varepsilon L}.
\end{align*}
The result follows by applying Lemma \ref{lemma:angle2d} after noticing that $|\vx-\vy|$ is the length of the shortest edge and that $|\vz-\vx|$ is the radius of the circumcirle. \jeffrey{y and z should be bold font, correct? They were not. Please confirm this change is okay.}\johannes{yes to all boldfont changes}
	\end{proof}
\end{lemma}

\begin{remark}
	Note that for $d>2$ the same result is true, if we assume the circumcenter of the $d$-simplex, of which $\Delta$ is a face, is contained in $\Omega$ instead of the circumcenter of $\Delta$ itself. The proof is identical.  
    \jeffrey{you were using both $n$ and $d$ for dimension. Previously, $n$ was the number of points in the sampling. I suggest $d$ for dimension, and $n$ for the sample set size (which changes to $N$ in the algorithms, fix that too). Go through check that the notation is now consistent.} \johannes{done}	
\end{remark}

\begin{remark}
	While this result allows us to control the quality of 2D-triangulations of maximal Poisson-disk samplings, it can also be used to gauge how close a given Poisson-disk sampling is to being maximal.
\end{remark}

The previous Lemma only gives us bounds on all triangles, if their circumcenters are contained in $\Omega$. 
The next two Lemmas will give sufficient conditions to guarantee exactly this as long as $\Omega$ is a polytope.
\begin{lemma}\label{lemma:circinside}

	Let $\Omega\subset\mathbb{R}^2$ be a polygonal region and $X$ a maximal Poisson-disk sampling containing all vertices of $\Omega$. Let the inhibition radius $r(\vx,\vy)$ be defined like in the previous lemma. 
	Furthermore, let the coverage radius of $X\cap \delta \Omega$ fulfill $R^\delta(\vx,\vy)<\frac{r(\vx,\vy)}{\sqrt{2}(1+L)}$, i.e., $|\vx-\vy|<{\frac{\sqrt{2}}{1+L}}r(\vx,\vy)$ for all $\vx,\vy\in \delta\Omega\cap X$. 
    Then the circumcenter of all triangles in the Delaunay triangulation of $X$ are contained in $\bar{\Omega}$.
	\begin{proof}
		Suppose this claim is wrong.
        Let $\Delta$ be a triangle in the Delaunay triangulation with circumcenter $\vz\notin \Omega$.\jeffrey{Changed z to a bold font. Please confirm ok.}
        For this to be possible the circumcircle needs to be cut into at least two pieces by $\delta\Omega$, separating $\vz$ and the vertices of $\Delta$.\jeffrey{Changed z to a bold font. Please confirm.} 
        Since $\Delta$ is part of a Delaunay triangulation and all vertices of $\Omega$ are part of the sampling, this is done by (at least) one segment of a straight line, i.e., $\delta \Omega$ contains a secant of the circumcircle.

		Next, let $\vb_1,\vb_2\in \delta\Omega$ be the two boundary points closest to the circumcircle on either side of that line segment and let $B$ be the disk bounded by the circumcircle. 
        Note that $\bar{B}\cap\Omega$ contains $\Delta$ and  is itself entirely contained in the disk of radius $\frac{1}{2}|\vb_1-\vb_2|<\frac{1}{2}{\frac{\sqrt{2}}{1+L}}r(\vb_1,\vb_2)$ around $\frac{1}{2}(\vb_1+\vb_2)$.

		Now, let $\vx\notin \{\vb_1,\vb_2\}$ be a vertex of $\Delta$ and let $\vb\in\{\vb_1,\vb_2\}$ be the point of the two, that is closer to $\vx$. 
        We already established that $\vx$ lies within the just mentioned ball around $\frac{1}{2}(\vb_1+\vb_2)$.
		Let $\vx_p$ be the projection of $\vx$ onto the line segment connecting $\vb_1$ and $\vb_2$. \jeffrey{Changed $b_1$ and $b_2$ to a bold font. Please confirm.}
        Then $|\vx-\vx_p|< \frac{1}{2}{\frac{\sqrt{2}}{1+L}}r(\vb_1,\vb_2)$, 
        because $\vx$ lies within the circle of that radius, $|\vx_p-\vb|< \frac{1}{2}{\frac{\sqrt{2}}{1+L}}r(\vb_1,\vb_2)$, because $\vb$ is the closer of the two points ${\vb_1,\vb_2}$ and therefore  \begin{align}
		\label{eq:rhob}
		|\vx-\vb|=\sqrt{|\vx-\vx_p|^2+|\vb-\vx_p|^2}< \frac{r(\vb_1,\vb_2)}{1+L}\le \frac{\rho(\vb)}{1+L}\le \rho(\vb).\end{align}
		Since $|\vx-\vb|\ge \min(\rho(\vx),\rho(\vb))$ this implies $|\vx-\vb|\ge \rho(\vx)$.
		However assuming this and applying the Lipschitz condition on \eqref{eq:rhob} gives us 
		\begin{align*}
			|\vx-\vb| < \frac{\rho(\vb)}{1+L}\le \frac{1}{1+L}\left(\rho(\vx)+L|\vx-\vb|\right)\le \frac{1+L}{1+L}|\vx-\vb|,
		\end{align*}
		which is a contradiction.
		\end{proof}
\end{lemma}
\begin{remark}
	Although Lemma \ref{lemma:circinside} generalizes to higher dimensions, it is not very practical because it is difficult to guarantee the bounds on $R^\delta$ if the boundary is more than 1-dimensional.
	However, it is still possible to get some bounds on the radii of the circumcircles and then, using Lemma \ref{lemma:angle2d}, on the angles, if the distance of non-boundary points is greater than some lower bound $a>0$. \jeffrey{$d$ was already used to denote spatial dimension, as well as using $n$ to denote for spatial dimension. I changed the $d \rightarrow a$  and $n \rightarrow d$. Please double check}\johannes{done}
    In fact, using notation from the previous proof, let $\Delta$ again be an $d$-simplex with circumcenter outside of $\Omega$ and $\vx \notin \delta\Omega$ on of its points.
    Since the circumsphere of any simplex in a Delaunay triangulation does not contain any other points the radius of the intersection with $\delta \Omega$ is bounded by $R^\delta$.
	Using simple geometric arguments, one can show that this feature forces the radius of $\Delta$'s circumsphere $R$ to fulfill the following inequality \begin{align} \label{eq:outsidecenter}
		R^2\le (R-a)^2+\left(R^\delta\right)^2 \Rightarrow R\le \frac{a^2+\left(R^\delta\right)^2}{2a}\le \frac{\left(R^\delta\right)^2}{a}.
	\end{align}  
	If $R^\delta(\vx,\vy)<r(\vx,\vy)$ there is a lower bound on $a$, continuously depending on $R^\delta$, solely due to the fact, that we have a Poisson-disk sampling. If $R^\delta=R^\delta(\vx_p,\vb)$ is any bigger, $a$ needs to be bounded artificially. This implies that the angle bounds change continuously, if the conditions for Lemma \ref{lemma:circinside} cannot be met they still can be relatively controlled by the choice of the artificial bound on $a$.
\end{remark}

Under the conditions of the previous Lemmas, the simplices of the two-dimensional Delaunay triangulation resulting from the maximal Poisson-disk sampling are guaranteed to only have well-behaved triangular faces. 
In three or more dimensions, however, this does not imply that the simplices themselves are well-behaved.
It is possible for a 3D Delaunay triangulation to contain slivers, which are tetrahedra whose four points are all positioned approximately on the equator of their circumsphere. 
In \cite{exudation} slivers are characterized as tetrahedra whose points are all close to a plane and whose orthogonal projection onto that plane is a quadrilateral. 
In \cite{bern} slivers are equivalently classified as tetrahedra with a dihedral angle close to $180^\circ$ containing their own circumcenter.
All faces of a sliver can be equilateral triangles, yet they can still have arbitrarily small dihedral angles.
This latter feature can cause large numerical errors in the computational physics simulations performed on the mesh.
Therefore, slivers should be eschewed as much as possible.
While slivers cannot be entirely avoided, one can show that if the points $\vx,\vy,\vz,\vw$ of a maximal Poisson-disk sampling form a sliver, the distance between $\vw$ and the plane spanned by $\vx,\vy,\vz$ is very small~\cite{exudation}.
We use this property to minimize slivers around the DFN and the faces of the surrounding volume by using the 2D sampling of the DFN as the seed of the 3D sampling and enforce a minimum distance between the DFN and subsequent points in the 3D sample.  
This procedure leads to slivers being scarce in the initial 3D triangulation.
Moreover, given any three points in the mesh, the probability a fourth point will produce is a sliver quite small in practice.   
Therefore, if the points producing slivers are removed and resampled, the resulting new triangulation will most likely have fewer slivers than the previous one. 
We tested this rejection/resampling algorithm and found it improves the overall quality of the mesh.


\section{Method: The near-Maximal Algorithm for Poisson-disk Sampling (nMAPS)}\label{sec:methods}

\subsection{nMAPS Workflow overview}

\begin{figure}
{\newpage
\hspace{-2cm}
\begin{minipage}[c]{0.45\textwidth}
    \vspace{-2cm}
    \begin{tikzpicture}
        \node (initial) at (.4,0) {\includegraphics[width=.2\textheight]{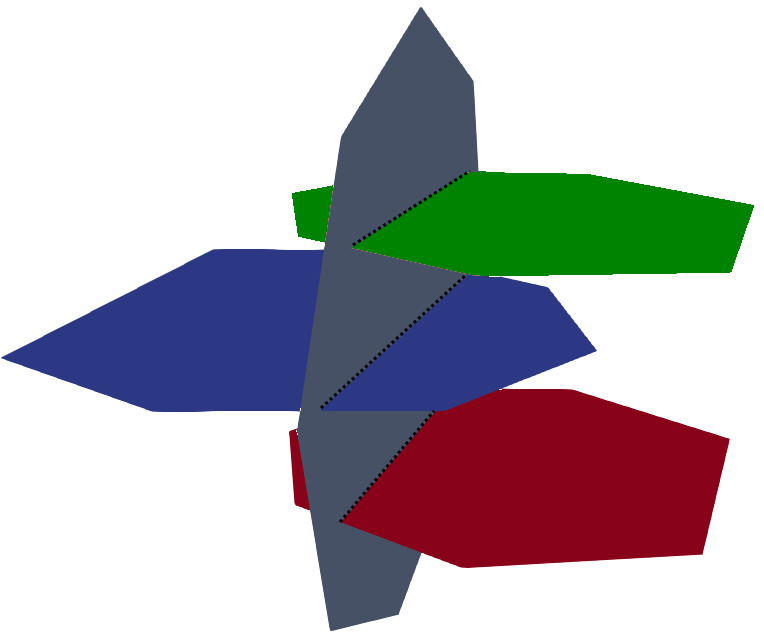}};
        \node (decomp) at (0, -4.5) {\includegraphics[width=.3\textheight]{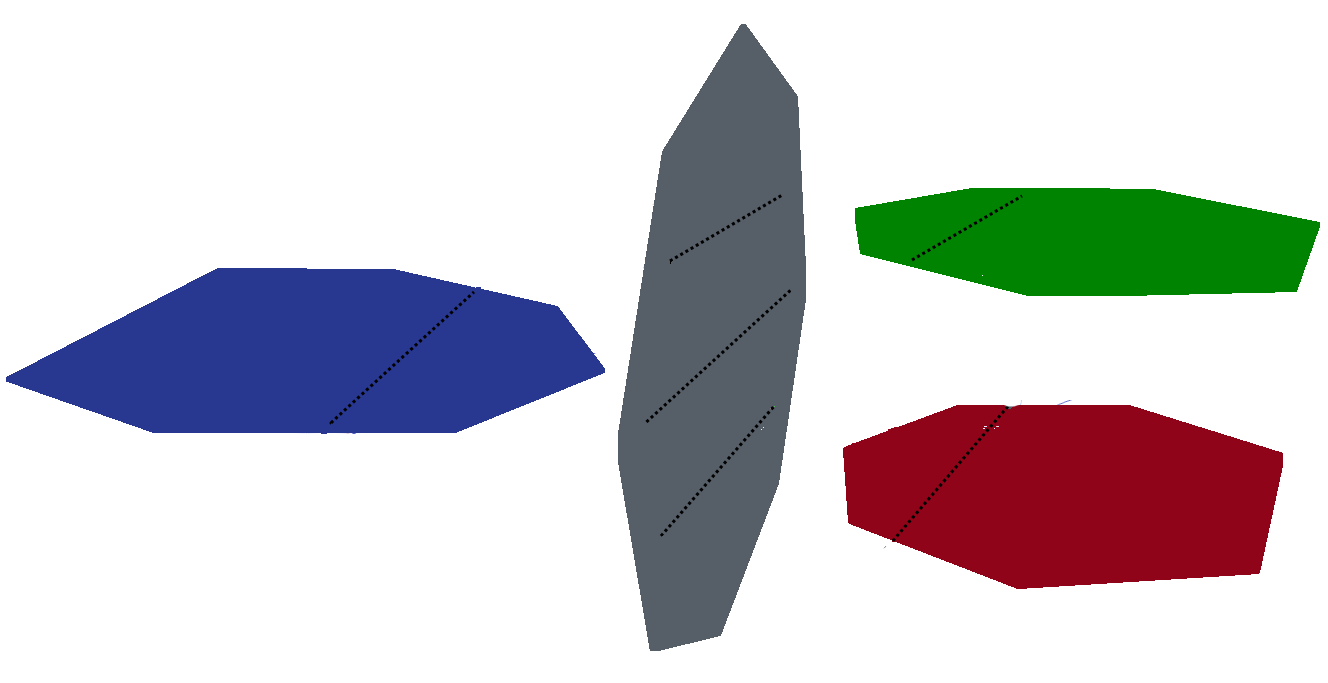}};
        \node (decomp) at (-3, -9) {\includegraphics[width=.3\textwidth]{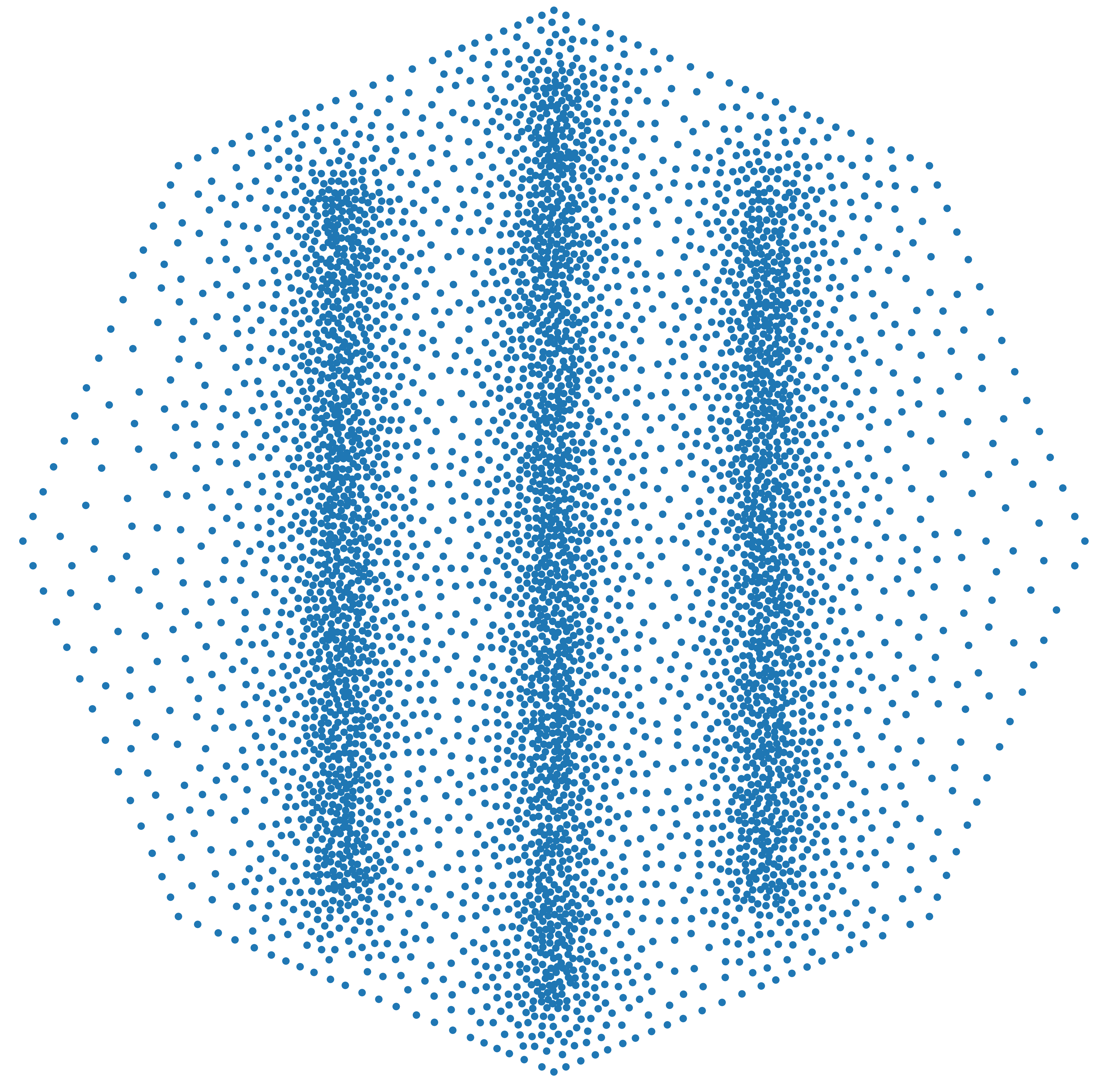}};
        \node (decomp) at (3, -9) {\includegraphics[width=.3\textwidth]{graphics/individual}};
        \node (decomp) at (0, -9) {\includegraphics[width=.3\textwidth]{graphics/individual}};
        \node (decomp) at (0.4, -13) {\includegraphics[width=.6\textwidth]{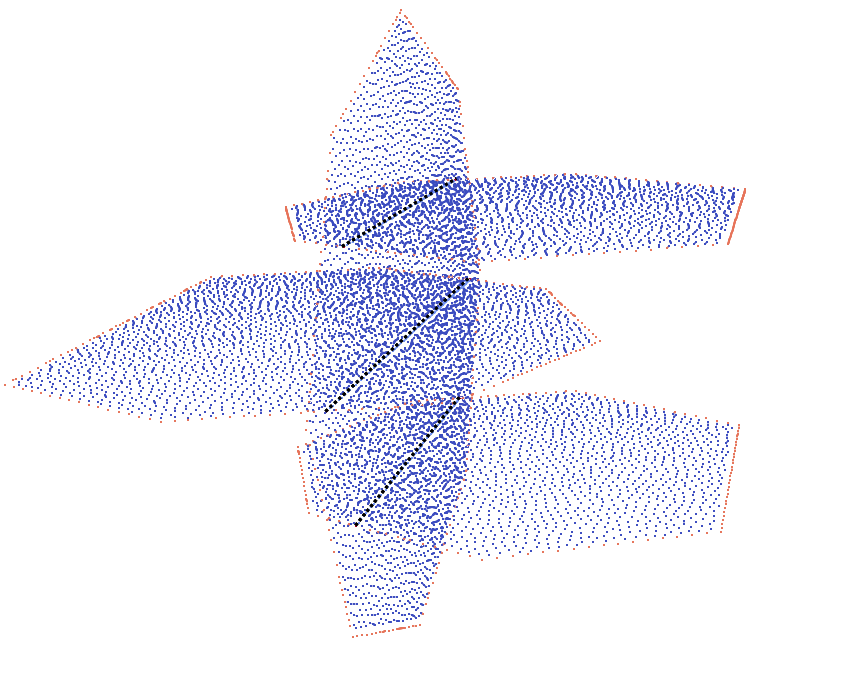}};
        \node (decomp) at (-1.2, -18) {\includegraphics[width=.7\textwidth]{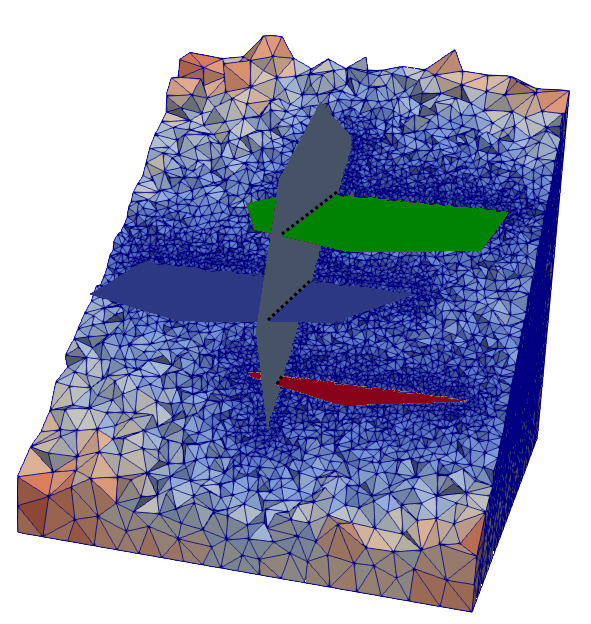}};

        \draw[->,black,ultra thick] (-1.,-2) -- (-1.5,-3)
        node[midway,fill=white]{\hspace{.5cm }};
        \draw[->,black,ultra thick] (1.4,-2) -- (1.9,-3)
        node[midway,fill=white]{        \hspace{.5cm }};
        \node[fill=white] () at (.4,-2.5) { 1. Deconstruct the DFN into individual fractures};

        \draw[->,black,ultra thick] (-2.2,-6.5) -- (-2.7,-7.5)
        node[midway,fill=white]{\hspace{.5cm }\vspace{1cm}};
        \draw[->,black,ultra thick] (-.0,-6.5) -- (-0,-7.5)
        node[midway,fill=white]{\hspace{.5cm }\vspace{1cm}};
        \draw[->,black,ultra thick] (2.,-6.5) -- (2.5,-7.5)
        node[midway,fill=white]{\hspace{.5cm }\vspace{1cm}};
        \node[fill=white,text width = 9cm] () at (0,-7) {2. Sample points each fracture independently (A-D)};

        \draw[<-,black,ultra thick] (-1.,-11.5) -- (-1.5,-10.5)
        node[midway,fill=white]{\hspace{.5cm }};
        \draw[<-,black,ultra thick] (1.4,-11.5) -- (1.9,-10.5)
        node[midway,fill=white]{        \hspace{.5cm }};
        \node[fill=white] () at (.4,-11.) {3. Re-assemble the DFN and eliminate conflicts};

        \draw[->,black,ultra thick] (0.,-15) -- (-0.3,-16)
        node[midway,fill=white]{\hspace{.5cm }};
        \node[fill=white] () at (-.2,-15.5) {4. Sample points in the surrounding volume};

        \node[text width=2.5cm,align=center] at (2,-17.5) {5. Recursively remove slivers};
        \end{tikzpicture}

\end{minipage} \hspace{3cm}
\begin{minipage}[c]{0.45\textwidth}
    \vspace{-4cm}
    
    \begin{tikzpicture}
    \node (init) at (5,0) {\includegraphics[width=0.7\linewidth]{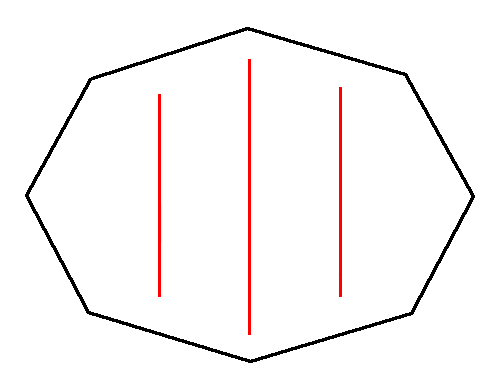}};
    \node (init) at (3,-4) {\includegraphics[width=0.7\linewidth]{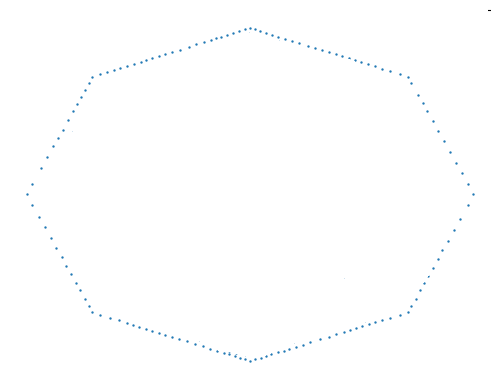}};
    \node (init) at (5,-8) {\includegraphics[width=0.7\linewidth]{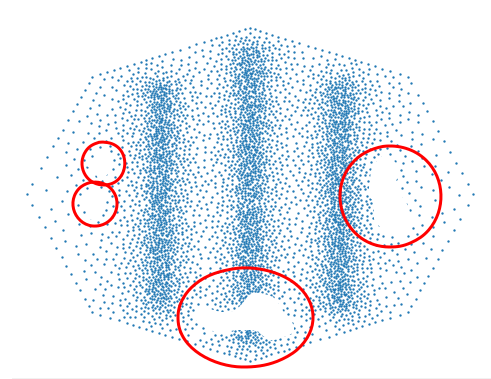}};
    \node (init) at (3,-12) {\includegraphics[width=0.7\linewidth]{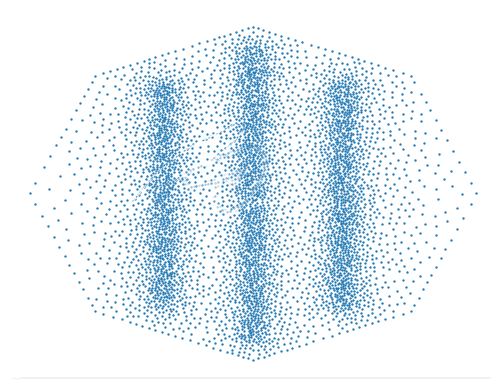}};
    
    \draw[->,black,ultra thick] (4,-1.5) -- (3.7,-2.5);
    \draw[->,black,ultra thick] (3.7,-5.5) -- (4,-6.5);
    \draw[->,black,ultra thick] (4,-9.5) -- (3.7,-10.5);
    
    \node[text width=5cm] () at (4,2) {A. Identify boundary and intersections};
    \node[fill=white] () at (3.85 ,-2.) {B. Sample points along boundary};
    \node[fill=white] () at (3.85 ,-6) {C. Use boundary points as the seed for nMAPS};
    \node[fill=white, text width = 3cm] () at (6.85 ,-11.5) {D. When the algorithm terminates, find the undersampled regions and sample new points therein. Repeat from step C.};    
    \end{tikzpicture}
    \vspace{0.5cm}

\caption{\label{fig:workflow} Overview of the nMAPS workflow from the creation of DFN to the final mesh (left 1-5) and during 2D-sampling (right A-D). } 
\end{minipage}}
\end{figure}

The nMAPS workflow contains the following high-level steps:
\begin{enumerate}
\item[1.] Generation of a DFN using dfnWorks \cite{dfn}. 
\item[2.] Deconstruct the DFN into individual fractures/polygons.
\item[3.] Generation of a 2D-variable-radii Poisson-disk sampling on each polygon.
\item[4.] Construct a conforming Delaunay triangulation of each fracture based on the Poisson-disk sample point set.
\item[5.] Merge individual fracture meshes into a single DFN. 
\item[6.] Generation of a 3D variable radii Poisson-disk sampling of the volume surrounding the meshed DFN.
\item[7.] Mesh generation of the resulting 3D sample set, identify low-quality tetrahedra and remove two of their nodes that are not located in the DFN mesh.
\item[8.] Repeat steps 6-7 with the remaining nodes as the seed set until all slivers are removed. 
\end{enumerate}
A graphical version of the workflow is presented in Figure \ref{fig:workflow}.
DFN generation was discussed in section~\ref{ssec:dfn}, so we begin with a description of the 2D sampling method.

\subsection{Two-Dimensional Sampling Method}
We generate a 2D Poisson-disk sampling in a successive manner using a rejection method.
First, we provide a general overview of the method.
This method can be performed on every fracture in the network independent of the other fractures. 
In each step, a new point candidate is generated.
It is accepted if it does not break the empty disk property with any of the already accepted points. 
For the sampling in two dimensions, we use a variable inhibition radius that increases linearly with distance from intersections on the fracture. 
Specifically, we reject a candidate point $\vy$, if there is an already excepted point $\vx$ such that the condition 
\begin{align}
|\vx-\vy|\ge r(\vx,\vy) = \min(\rho(\vx),\rho(\vy)) \label{eq:empty_disk}
\end{align}
is violated.
In this equation, $\rho(\vx)$ as a piecewise linear function given by
\begin{align}
\rho(\vx)=\rho(D(\vx))=\left\{\begin{matrix}
\frac{H}{2} &\text{for} & D(\vx)\le FH\\
A(D(\vx)-FH)+\frac{H}{2} &\text{for}& FH\le D(\vx) \le (R+F)H \\
 (AR+\frac{1}{2})H & & \text{otherwise}
\end{matrix}\right .\label{eq:rho2d}
\end{align}
Here $D(\vx)$ is the Euclidean distance between $\vx$ and the closest intersection. 
$H, A, R$, and $F$ \jeffrey{the parameters are in math mode here, but elsewhere they are not. I fixed what I found, but please go through and have another look}\johannes{looks like you found them all} are parameters that determine the global minimum distance between two points ($H/2$), the distance around an intersection where the local inhibition radius remains at its minimum value ($FH$), the global maximum inhibition radius ($ARH+H/2$), and the slope at which the inhibition radius grow with $D(\vx)$ ($A$).
Since $\rho(D)$ is piecewise linear, it is a Lipschitz-function with Lipschitz-constant $A$.

\begin{figure}
	\centering
	\includegraphics[height=0.4\linewidth]{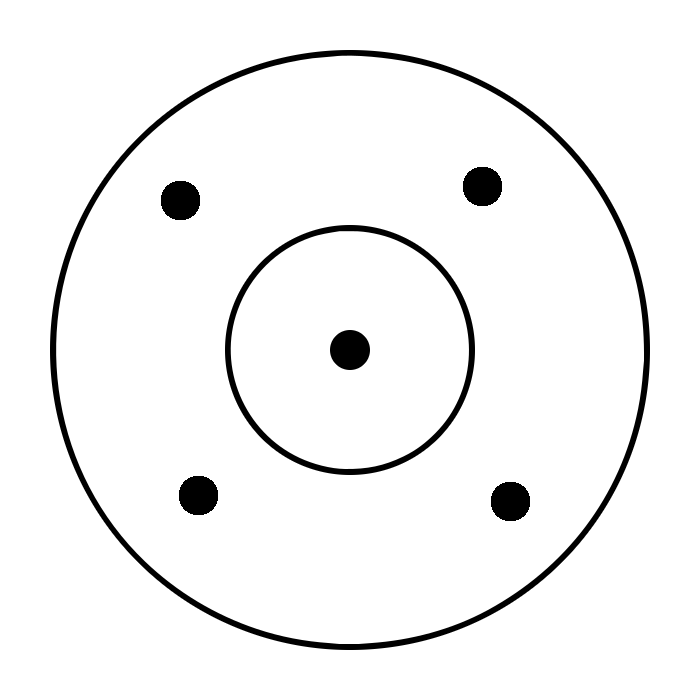}
	\caption{Visualisation of a single sampling step. The current point is in the center, the new candidates in the annulus ($k=4$). The inner circle is bounded by the inhibition radius of the current point. The outer circle is bounded by the maximum distance a point can be away from the current point if the Poisson-disk sampling was maximal. }
    \label{fig:annulus}
\end{figure}

If the sampling has a coverage radius $R(\vx,\vy)\le (1+\varepsilon)r(\vx,\vy)$ for some $\varepsilon>0$, then the conditions of \eqref{lemma:2} hold. 
To satisfy the conditions of \eqref{lemma:circinside} as well, and thereby ensure angle bounds on all triangles in a Delaunay triangulation, we begin by sampling points along the boundary and enforce a maximum distance of $\frac{r(\vx,\vy)}{\sqrt{2}(1+L)}$ between points.
Next, we generate new candidates by randomly sampling within an annulus around an already accepted point; this is illustrated in Figure \ref{fig:annulus}.  
The minimum distance another point could be from the center point and still preserve the empty disk property determines the inner radius.  
The maximum distance a point be from to the center in a maximal sampling determines the outer radius.
For our choice of inhibition radius, and assuming it has the same radius as the coverage radius, these distances are $\frac{\rho(\vx)}{1+A}$ for the inner radius and $\frac{2\rho(\vx)}{1-A}$ for the outer radius.

\begin{figure}
	\centering
(a)	\includegraphics[height=0.4\linewidth]{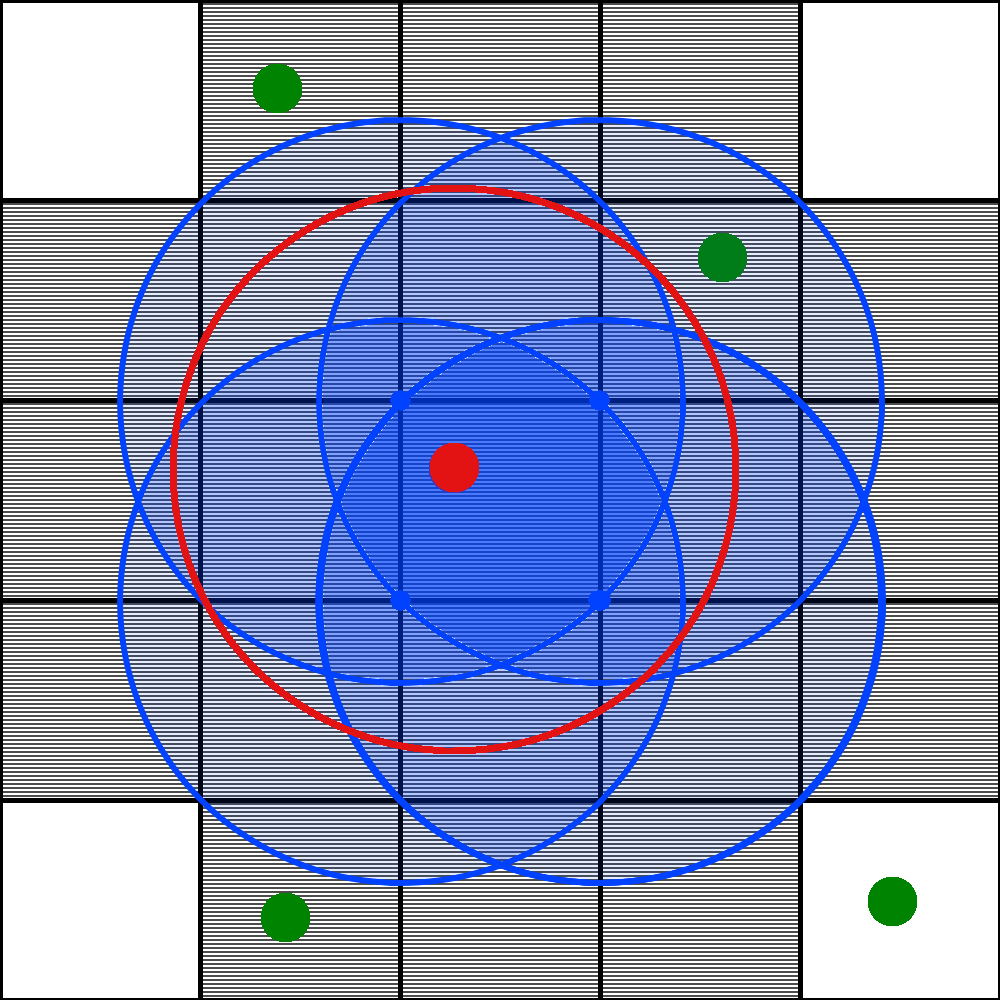} \hspace{0.5cm}
(b)		\includegraphics[height=0.4\linewidth]{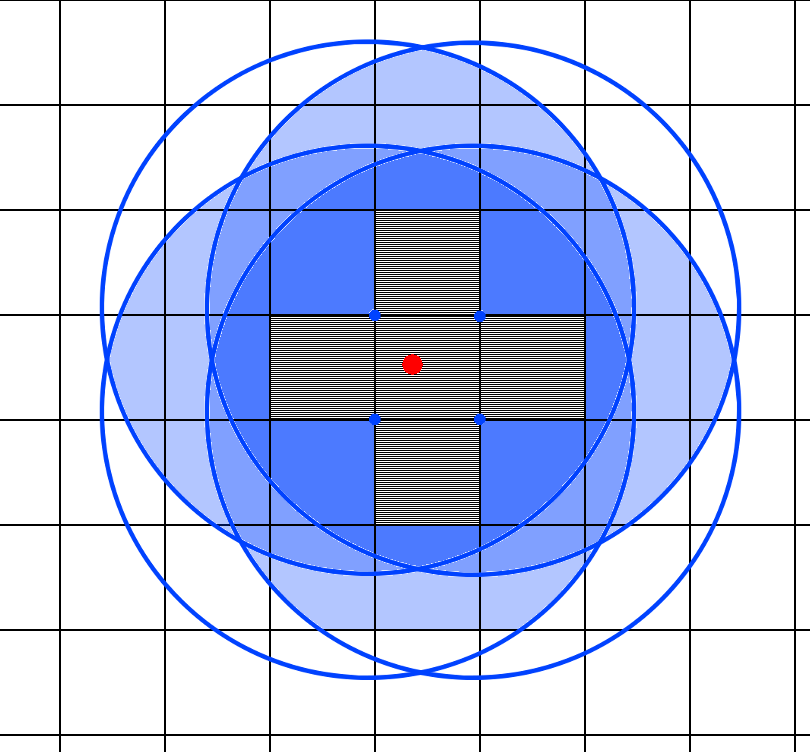}
	\caption{Visualization of how the grid is used to find possibly conflicting points. The new candidate is red, already accepted points are green, and cells containing conflicting points are shaded grey. The red circle shows the inhibition radius of the candidate, and blue circles show the furthest cells a point in the center cell could conflict with.}
	\label{fig:grid}
\end{figure}

We will now go over the individual steps of the 2D algorithm in detail. 
These steps are also presented in the pseudocode in Algorithm 1 (Section \ref{sec:pseudo}) and are illustrated in A-D of Fig.~\ref{fig:workflow}. 
The notation used in the pseudocode is found in the table at the start of the same section. 
First, a 1D Poisson-disk sampling along the boundary of the polygon is generated as a seed set (line \ref{line:seed}). 
Next, $k$ candidate points at a time (line \ref{line:k_samp}) around each already accepted point are sampled and determined whether they are accepted or not (lines \ref{line:rej} through \ref{line:reject_end}). 
$k$ is a positive integer and a user-defined parameter. 
We use cell lists to find points around a candidate that could potentially cause a violation of the empty disk. 
A visual depiction of these cells is provided in Fig.~\ref{fig:grid}(a). 
The size of these cells is chosen to contain at most one point, which allows us to skip distance calculations with points beyond a certain cutoff (line \ref{line:yinNp}). 
Unlike the previously mentioned algorithms, we label cells containing points as occupied \textit{and} those cells that are too close to an accepted point to contain the candidate point.  
Specifically, if a candidate $\vx$ lies in a cell $C$ and any other cell $D$ with $diam(C\cup D)\le r_{in}(\vx)$ is occupied, then $\vx$ is immediately rejected as it conflicts with the point in $D$ (line \ref{line:rej}). 
On the other hand, if $dist(C, D)>\rho(\vx)$, then there is no need to calculate the distance between $\vx$ and any potential element of $D$ because they cannot violate the empty disk property. 
An example of this property is shown in Fig.~\ref{fig:grid}(b).
We use this technique to our advantage in two ways.
First, it allows us to reject many candidates without calculating any distances to nearby points, which provides substantial speedup compared to previous methods and especially for large values of $k$.
Second, unmarked cells that contain space for another point are easily identified, which allows us to find under-sampled regions after the initial sweep terminates (line \ref{line:emptycell}). 
If a point is accepted, then it is added to the end of the sample set (line \ref{line:accept_start}).
Note that this newly accepted point will be used as a seed for sampling another $k$ point when its turn comes up in the queue. 
Thus, the method grows the sample set within the primary while loop.
If all $k$ candidates around a point are rejected, then we move on to the next already accepted point (line \ref{line:next_node}). 
The algorithm terminates once every accepted point has been used as a sampling center (line \ref{line:end}).
Next, we detect unmarked or under-sampled cells and randomly place points within them (line \ref{line:rand_new}). 
Then, the main algorithm is restarted with the seed set, including these newly added points, and ends once the termination criterion is met once again (line \ref{line:rerun}). 
While this process of resampling sweeps in under-sampled cells can be performed multiple times, we found that a single resampling was sufficient to increase the quality of the sampling to acceptable levels, and additional sweeps were unnecessary.

Once the point distribution is created, the conforming Delaunay triangulation method of Murphy et al.,~\cite{murphy2001point} is used to generate the mesh.
We recount the general idea of the method here for completeness.
To create a conforming Delaunay triangulation that preserves the lines of fracture intersections as a set of triangle edges to be created, it is sufficient that the circumscribed circle of each segment of the discretized line of intersection be empty of any other point in the distribution prior to connecting the mesh.
To achieve this sufficient condition, any point within the circumscribed circle of each segment of the discretized lines of intersection is removed from the point distribution. 
Next, a two-dimensional unconstrained Delaunay triangulation algorithm is used to connect the resulting point set. 
Because of the construction method, i.e., empty regions around the lines of intersection, the line segments that represent lines of fracture intersection naturally emerge in the triangulation. 
In turn, the Delaunay triangulation will conform to all of the fracture intersection line segments. 
Once every fracture polygon is triangulated, they are all joined together to form the mesh of the entire DFN.

\subsection{Three-dimensional Sampling Method}
nMAPS uses a similar method to that presented in the two-dimensional section to generate a point distribution in three-dimensional space. 
The primary difference is that candidates are generated on a spherical shell around accepted nodes instead of an annulus. 
The 3D variant of $\rho(\vx)$  given by
\begin{align}
\rho(\vx) = \rho(D(\vx))=\left
\{\begin{matrix}
\rho_{2}(\vx_p)  &\text{ for }&  D(\vx)\le F\rho_2(\vx_p)\\\\
A(D(\vx)-Fr_{2D}(\vx_p))+\frac{H}{2} & \text{ for }& F\rho_{2D}(\vx_p)\le D(\vx) \le \frac{\rho_{max}-r_{2}(\vx_p)}{A} \\\\
 \rho_{max}& & \text{otherwise} 
\end{matrix}
\right. \label{eq:rho3d}
\end{align}
where $\vx_p$ the fracture point closest to $\vx$, $\rho_2(\vx_p)$ is its 2D inhibition radius on the fracture, and $D(\vx)$ is the distance between $\vx$ and $\vx_p$.
Just as in the 2D case, equation \eqref{eq:rho3d} is a piecewise linear function in $D(\vx)$. 
It is constant within a distance of $\rho_2(\vx_p)F$ of the fracture network and increases linearly with a slope of $A$ until a given maximal inhibition radius of $\rho_{max}$ is reached. 
\jeffrey{How is $\rho_{max}$ defined? }\johannes{$\rho_{max} is a parameter$, in practice same as max radius in 2D, i.e. (AR+0.5)H. Shifting to an upper bound on the radius rather than the range over which we increase, because the inhibition radius on the fracture is not constant.} In practice these parameters can be chosen to be the same as their 2D counterparts.
One difference from the 2D case is that in addition to rejecting all candidates that violate the empty-disk property \eqref{eq:empty_disk}, we also reject a candidate $\vx$ if it is within a distance of $\rho(\vx)/2$ to a boundary or fracture. 
This latter piece prevents slivers from having three nodes located on a single fracture or the domain's boundary.
In turn, this limits the circumradius of tetrahedra with circumcenter outside of the domain; cf. lemma \ref{lemma:circinside} and subsequent remark for additional details.

The pseudocode of the 3D-sampling method provided is in Algorithm 2 in Section \ref{sec:pseudo}. 
The initial sampling process is identical to the 2D version, therefore, we will focus on the initialization and resampling, which is where the methods differ.
The seed set is made up of points sampled along the boundary of the 3D volume and those of the DFN generated by the 2D algorithm (line \ref{line:seed3d}).
Neighbor cells can still be used in the same way as in 2D to speed up the rejection of candidates. 
Unlike in 2D, a maximal Poisson-disk sampling in 3D does not guarantee sliver-free triangulation, which is why we do not use the cell lists to find undersampled cells in 3D. 
Instead, once the algorithm terminates, the resulting sampling is triangulated (line \ref{line:triangulate}), slivers identified (line \ref{line:sliver}), and two nodes of every sliver, with a preference for nodes that are neither on a boundary or a fracture, are removed (line \ref{line:rem_nodes}). 
While the definition of a sliver given earlier in Section \ref{sec:mpds} allows for a bit of leeway in what is considered a small or large dihedral angle, we successfully replaced tetrahedra with dihedral angles outside of $[8^\circ,170^\circ]$ and aspect ratios bigger than $0.2$.
It is worth noting that more traditional ways of sliver-removal like perturbation \cite{Perturbation} or exudation \cite{exudation} can break the empty disk property and are not used in nMAPS.
Then the algorithm is restarted with the remaining nodes as the seed set (line \ref{line:rerun3d}). 
This process is repeated until all slivers are removed (line \ref{line:term3d}).
With this approach, we have obtained triangulations with no elements of dihedral angles of less than  $8^\circ$, examples provided in the next section.
The method for generating the conforming mesh is similar to that for the 2D case, but spheres around triangle cells of the fracture planes are excavated. 


%

\subsection{Pseudocode for the 2D and 3D sampling algorithms\label{sec:pseudo}}
\begin{tabularx}{\textwidth}{l X}
\multicolumn{2}{l}{ \textbf{Notation for Pseudocodes:}}\\ 	\multicolumn{2}{l}{ \hrulefill}\\
	\multicolumn{2}{l}{ \textbf{Input:}}   \\
	\textbullet $D^3\subset\mathbb{R}^3$: & Cubical Domain ($\ast$) \\   
	\textbullet $DFN\subset D^3$: & Generated by {\sc dfnWorks} ($\ast$) \\
	\textbullet	$F_l \subset \mathbb{R}^3$  :& $l$-th fracture of the DFN  \\
    \textbullet  $q_{l,m}^{(1)}$ and $q_{l,m}^{(1)}$:& Endpoints of intersection between fractures $F_l$ and $F_m$ \jeffrey{is $q_{l,m}$ used? I can't find it elsewhere. Remove if not.} \johannes{It's implicitly used to calculate r(x,y). If someone used a different r(x,y) they might not need it, or need sth. else. The same is true for $H,F,R,A$ though. Remove? } \\  
	\multicolumn{2}{l}{\textbf{User defined parameters: }}   \\
\textbullet  $H/2$: & Minimum distance between points \\
 \textbullet  $F$ : & $HF$ distance round intersections with constant density\\
 \textbullet  $R$:& $ARH+H/2$ is maximum distance between points  \\
 \textbullet $A$:& Maximum slope of inhibition radius  \\
  \textbullet$k$:& Number of concurrently sampled candidates\\
 \multicolumn{2}{l}{\textbf{Additional notation:}}   \\
 \textbullet $G$:& Square cells covering $F_l$ with  $diam(g)\le H/2$ for all $g\in G$.\\
 \textbullet $\rho(\vx)$: & Piecewise linear function defined in \eqref{eq:rho2d}-2D or \eqref{eq:rho3d}-3D \\
 \textbullet $r(\vx,\vy)$: & Inhibition radius $\min(\rho(\vx),\rho(\vy))$ \\
 \textbullet $R(\vx,\vy)$: & Coverage radius \\
\textbullet$C(\vx) \in G$:& Grid cell containing the point $\vx$.\\
 \textbullet $N^+(\vx)$ :& $\{g\in G: dist(C(\vx),g)\le \rho(\vx)\}$: Cells that can contain points $\vy$ with $|\vx-\vy|\le r(\vx,\vy)$ \jeffrey{Changed $y$ to $\vy$. Please confirm this is okay}\\
 \textbullet $N^-(\vx)$: &$\{g\in G: diam(g\cup C(\vx))\le \frac{\rho(\vx)}{1+A} \}$: Cells, where for all their points $y$ $|\vx-\vy|\le r(\vx,\vy)$\\
 \textbullet $G_{occ}$& $\bigcup\limits_{\vx\in X} N^-(\vx)$: Cells on which $X$ is already maximal.\\
 \textbullet $\mathcal{T}(X)$: & Delaunay triangulation of $X$ ($\ast$)\\
 \textbf{Output:} &  \\
 \textbullet  $X$ :& Poisson-disk sampling on the $l$-th fracture \\ 
 	\multicolumn{2}{l}{\hrulefill}\\
 	$(\ast)$: Only used in 3D sampling 
\end{tabularx}

\begin{algorithm}[H]	
	\caption{1}{\textbf{2D Poisson-disk sampling}\label{alg:2d}}

	\begin{algorithmic}[1]		
		\STATE\textbf{Initializing:}
		\STATE $X \subset F_l$ \hfill $\triangleright$ Generate a 1D Poisson-disk sampling with $R(\vx,y)\le \frac{r(\vx,\vy)}{\sqrt{2}(1+L)}$\\ \hfill along boundary $\delta F_l$ as seed. \label{line:seed}
		\FOR{$\vx\in X$}
		\STATE $G_{occ}\leftarrow G_{occ}\cup N^-(\vx)$ \hfill $\triangleright$ Initialize occupied cells
		\ENDFOR
		 
		 \hspace{0cm}
		\STATE\textbf{Sampling:}
		\STATE $i\leftarrow 1$ \hfill $\triangleright$ Start sampling at first accepted point
		\STATE $n\leftarrow |X|$ \hfill $\triangleright$ Will increase as more points are accepted
		\WHILE{$i\le n$}  \label{line:sampling}
		
		\REPEAT
		\FOR{$j \in \{1,...,k\}$} \label{line:gen_candidates}
		\STATE	$\vp_j \in F_l$ \hfill $\triangleright$ Generate $k$  new candidate points on the annulus around $\vx_i$\label{line:k_samp}
		\IF{$C(\vp_j) \in G_{occ}$} \label{line:rej}
		\STATE reject $\vp_j$ 	\hfill $\triangleright$ Cell already blocked by existing point's inhibition radius
		\ELSE 
		\FOR{$\vy\in N^+(\vp_j)$}\label{line:yinNp}
		\IF{$|\vp_j-\vy|<r(\vp_j,\vy)$}
		\STATE reject $\vp_j$ 	\hfill $\triangleright$ Empty disk property violated
		\ENDIF
		\ENDFOR
		\ENDIF \label{line:reject_end}
		\IF{$p_j$ was not rejected} \label{line:accept_start}
		\STATE  $X \leftarrow X\cup \{\vp_j\}$\hfill $\triangleright$ Accept $\vp_j$ and add it to the sampling set
		\STATE $G_{occ}\leftarrow G_{occ} \cup N^-(\vp_j)$ \hfill $\triangleright$ Update occupied cells	
		\STATE $n\leftarrow n+1$ 	\hfill $\triangleright$ Ensures sampling around newly accepted \label{line:accept_end} points				
		\ENDIF		
		\ENDFOR 
		\UNTIL{All $k$ of the $\vp_j$ are rejected}\label{line:all_k_rejected}
		\STATE $i\leftarrow i+1$ \hfill $\triangleright$  Start sampling around next accepted point \label{line:next_node}
		\ENDWHILE \hfill $\triangleright$ Terminate here or start resampling \label{line:end}
	\end{algorithmic}
\end{algorithm}

\begin{algorithm}[H]
\captionsetup{labelformat=empty}
\caption{}{\textbf{Continuation of Algorithm 1 (2D Resampling)}}

	\begin{algorithmic}[1]	
		\setalglineno{31}	
		\STATE  \textbf{Resampling:} (Optional)
		
		\FOR{$C\in G \setminus G_{occ}$}\label{line:emptycell}
		\STATE $\vp \in C$  \hfill $\triangleright$ Generate a random candidate on each cell \label{line:rand_new}
		\FOR{$\vy\in N^+(\vp)$}
		\IF{$|\vp-\vy|<r(\vp,\vy)$}	
		\STATE reject $\vp$ \hfill $\triangleright$ Empty disk property violated
		\ENDIF
		\ENDFOR		
		\IF{$\vp$ was not rejected} 
		\STATE  $X \leftarrow X\cup \{\vp\}$\hfill $\triangleright$ Accept $\vp$ and add it to the sampling set
		\STATE $G_{occ}\leftarrow G_{occ} \cup N^-(\vp)$ \hfill $\triangleright$ Update occupied cells 	
		\STATE $n\leftarrow n+1$

		\ENDIF

		\ENDFOR
		\STATE
		\STATE \textbf{Rerun algorithm from line \ref{line:sampling}} (Note: $i$ is not reset)\label{line:rerun}
	\end{algorithmic}

\end{algorithm}

\begin{algorithm}[H]	
	\captionsetup[FLOAT_TYPE]{labelformat=simple}
	\caption{2}{\textbf{3D Poisson-disk sampling + resampling}\label{alg:3d}}

	\begin{algorithmic}[1]
		\STATE \textbf{Initializing:}
		\STATE $X\subset D^3$ \hfill $\triangleright$ Use Algorithm 1 to generate a Poisson-disk sampling on $\delta D^3$ and the DFN \label{line:seed3d} \label{line:init3d}
		\FOR{$\vx\in X$}
		\STATE $G_{occ}\leftarrow G_{occ}\cup N^-(\vx)$ \hfill $\triangleright$ Initialize occupied cells \label{line:hfhf}
		\ENDFOR

		\STATE	 \textbf{Sampling:} \label{line:sampling3d}
		\STATE The sampling process in 3D is the same as the 2D method (Algorithm 1) except for the following: 
		\begin{itemize}
			\item[A.] New candidates are generated on a spherical shell instead of an annulus
			\item[B.] A candidate $p\notin \delta D^3$ is rejected if $dist(p,\delta D^3)<\rho(p)/2$
		\end{itemize}
		\STATE \textbf{Resampling:}  (Optional)
		
		\REPEAT 
		\FOR{$T\in \mathcal{T}(X)$} \label{line:triangulate}
		\IF{$T$ is a sliver} \label{line:sliver}
		\STATE $X\leftarrow X\setminus\{\vx,\vy\}$, where $\vx,\vy\in T$ are two random points not contained in the boundary or the DFN \label{line:rem_nodes}
		 \hfill $\triangleright$ Minimum distance of points to DFN and boundary (7B above) ensures that this is possible.
		\ENDIF
		\ENDFOR

		\STATE \textbf{Rerun algorithm from line \ref{line:sampling3d}}\label{line:rerun3d}
		\UNTIL{$\mathcal{T}(X)$ contains no more slivers.}\label{line:term3d}
		
	\end{algorithmic}

\end{algorithm}

\section{Results}\label{sec:results}

In this section, we provide examples of nMAPS and compare its performance with previous methods.

\subsection{Two-dimensional Examples}

\begin{figure}
	\centering
            \includegraphics[width=1\linewidth]{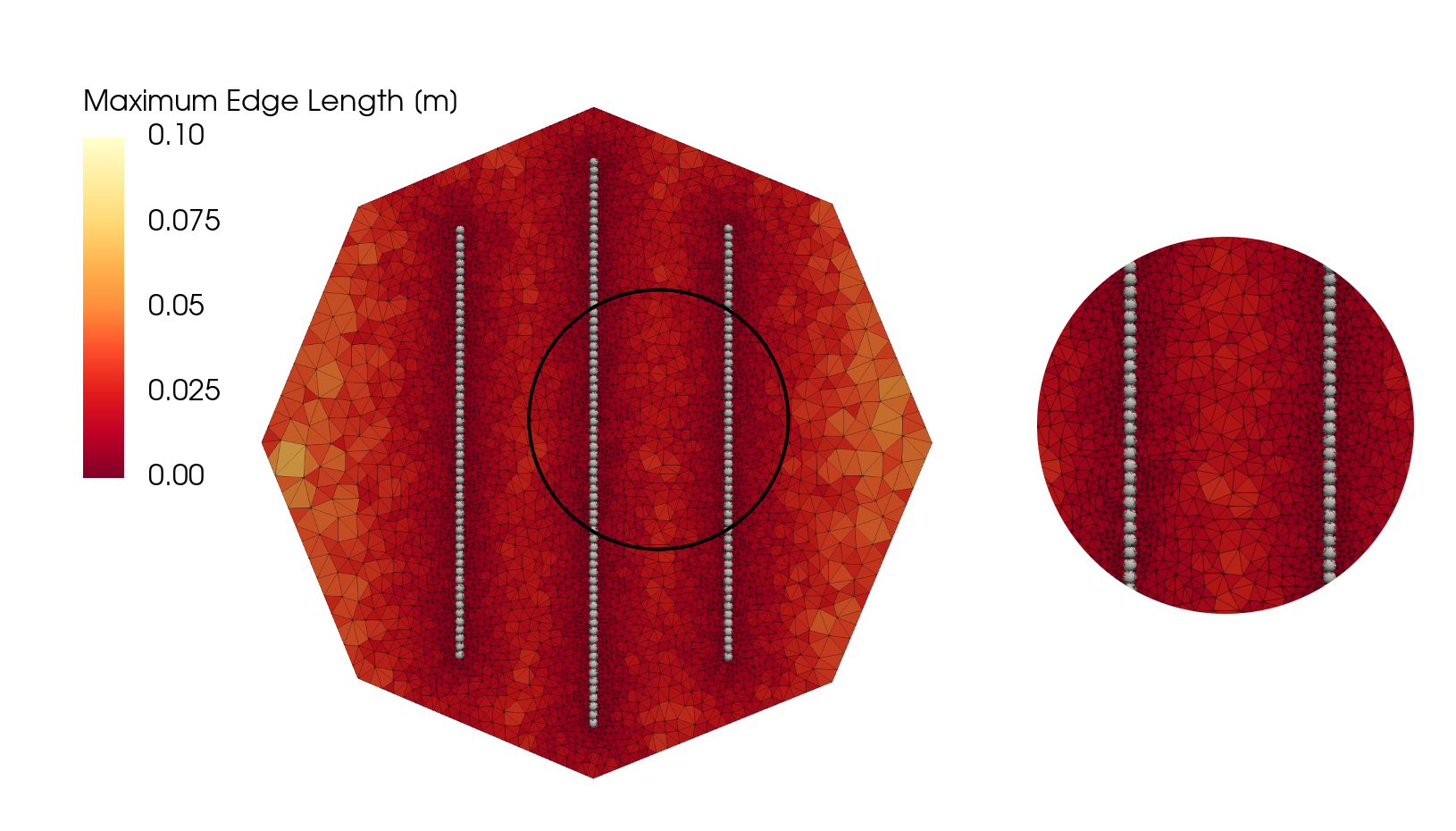}
	\caption{Triangulation of variable radii Poisson-disk sampling on fracture with three intersections. Parameters $H=0.01$, $R=40$, $A=0.1$, and $F=1$. Triangles are colored according to their maximum edge length. The lines of intersection are shown as spheres. }
	\label{fig:varedgmaxfin}
\end{figure}

\begin{figure}
	\centering
	\includegraphics[height=0.4\linewidth]{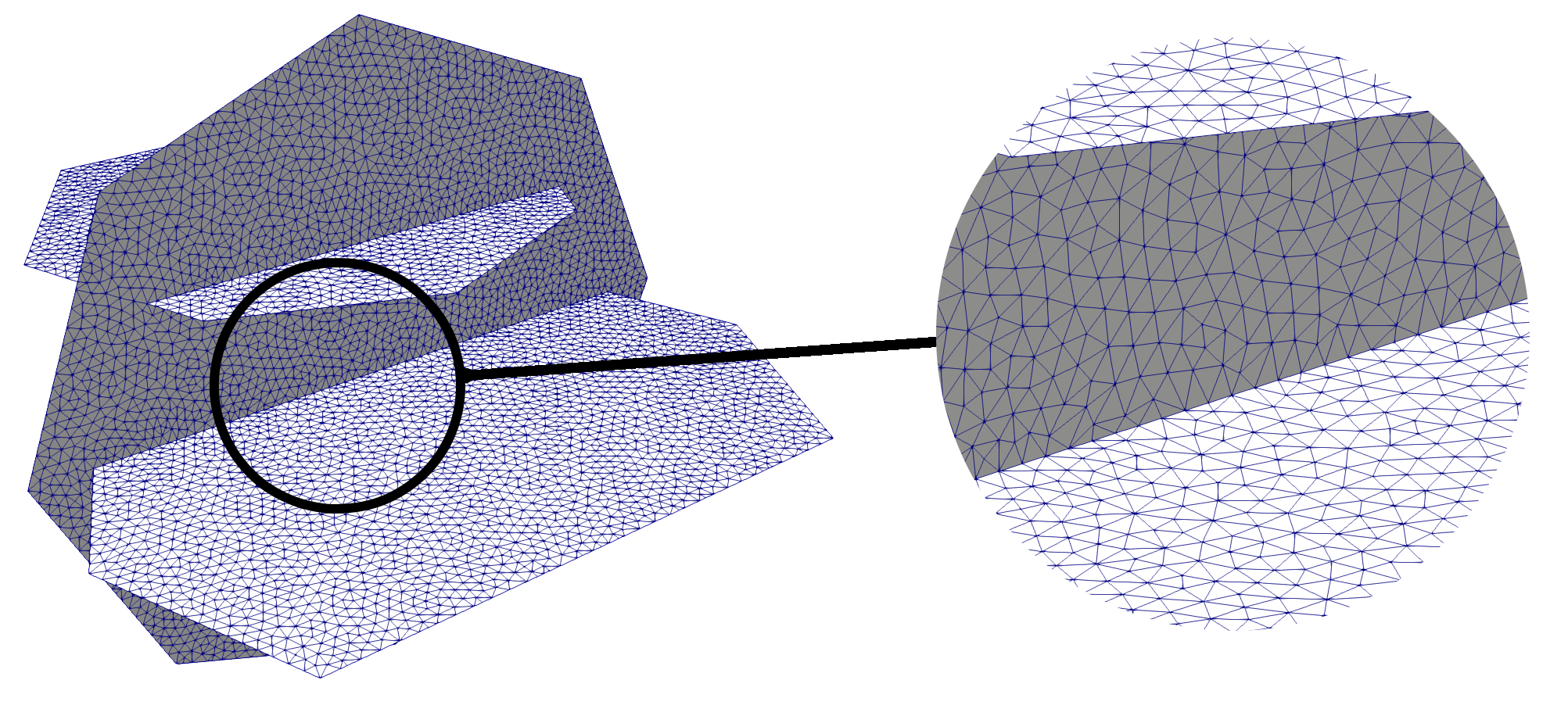}
	\caption{Triangulation of a uniform resolution Poisson-disk sampling based mesh reassembled into the whole DFN.}
	\label{fig:2d-together}
\end{figure}

We begin with a network composed of four disc-shaped fractures. 
Figure \ref{fig:varedgmaxfin} shows the triangulation produced using a variable radius sampling on a single fracture that contains three intersections. 
Triangles are colored by their maximum edge length to demonstrate how their size increases with distance from the intersections.  
Figure \ref{fig:2d-together}, shows the triangulation of a constant-radius sampling (uniform resolution mesh) on the fracture but re-assembled into the network.

\begin{figure}
	\centering
	\includegraphics[height=0.4\linewidth]{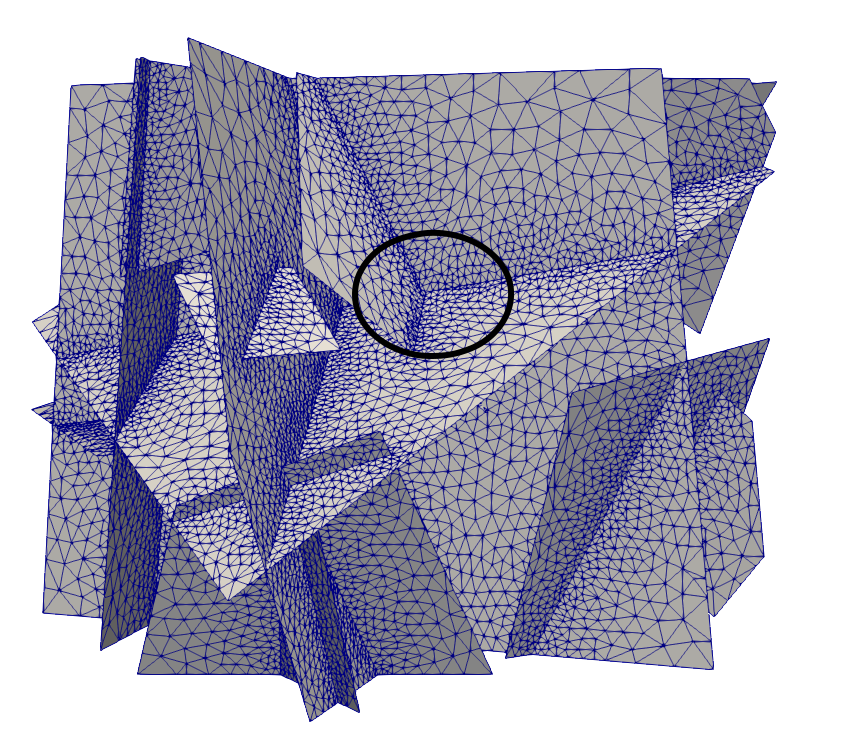} 
	\includegraphics[height=0.3\linewidth]{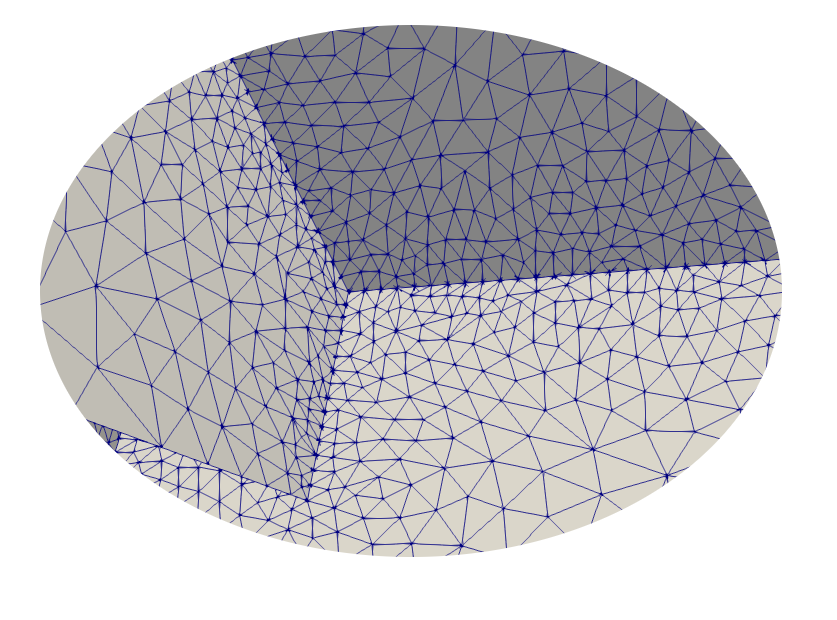}
	\caption{Triangulation of a variable radii Poisson-disk sampling reassembled into original DFN. The network contains 25 fractures whose radii are generated from an exponential distribution with a decay exponent of 0.3. Parameters used in the sampled: $H=0.1$,$R=40$,$A=0.1$,$F=1$. The mesh contains 23195 nodes and 47367 triangles. The minimum angle is $\ge 25^\circ$, the maximum angle $\le 120^\circ$, and all aspect ratios are $\ge0.47$. }
	\label{fig:togetherzoomexp}
\end{figure}

\begin{figure}
	\centering
	\begin{tikzpicture}
	\node (A) at (0,0) {\input{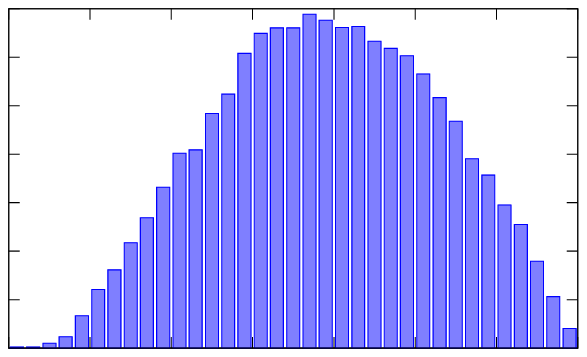}};
	\node (B) at (0,-5.5) {\input{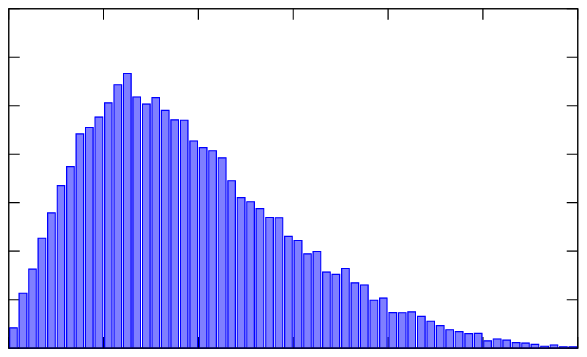}};
	\node (B) at (0,-11) {\input{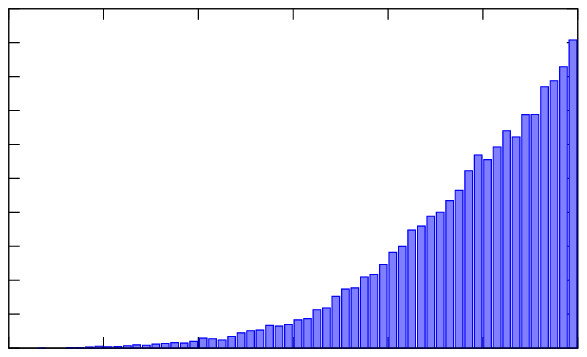}};
	\node (a) at (-4.,2.3) {(a)};
	\node (a) at (-4.,-3.2) {(b)};	
	\node (a) at (-4.,-8.7) {(c)};		
	\end{tikzpicture}
	\caption{Histograms of selected quality measures of the triangulation of variable radii Poisson-disk sampling on a fracture with three intersections. Parameters: $H=0.01$, $R=40$, $A=0.1$, and $F=1$. (a): Minimum angle ($\ge 25^\circ$), (b): Maximum angle ($\le 120^\circ$), (c): Aspect ratio ($\ge 0.47$)}
	\label{fig:variablemaxangle}
\end{figure}

Figure \ref{fig:togetherzoomexp} shows a meshed network that contains 25 fractures whose radii are generated from an exponential distribution with a decay exponent of 0.3.
The largest number of intersections on a fracture is eight within the network; this is not a constraint of generation nor the sampling technique.
The inhibition radius parameters are $H=0.1$, $A=0.1$, $F=1$, and $R=40$.
The mesh contains 23,195 nodes and 47,367 triangles.
The quality of the triangulation is presented in the histograms shown in Figure \ref{fig:variablemaxangle}: distribution of (a) minimum angle, (b) maximum angle, and (c) aspect ratio. 
With the exception of two elements, all of the triangles have a minimum angle greater than $27^\circ$.
The theoretical minimum angle in a maximal Poisson-disk sampling with Lipschitz constant $A=0.1$ is $27.04^\circ$. 
The two exceptions are $25^\circ$ and $26^\circ$.
In terms of the maximum angle, there are very few triangles with angles larger than $110^\circ$ and none larger than $120^\circ$. 
The largest maximum angle theoretically possible in a maximal Poisson-disk sampling with this Lipschitz-constant is $125.92^\circ$. 
The vast majority of aspect ratios are greater than $0.8$ with only a marginal number of triangles having a value less than $0.6$ and none less than $0.47$.

\subsection{Run Time Analysis}
\begin{figure}
	\centering
	\input{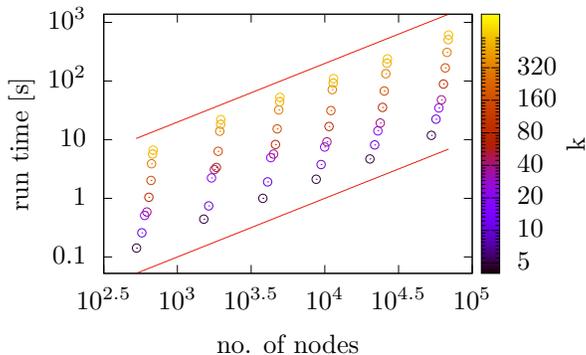}
	\caption{Run time of nMAPS plotted as a function of points sampled prior to the resampling process. 
    Data points are generated using the same DFN. 
    Different point densities are generated by changing the minimum inhibition radius  $\frac{H}{2}$ between every pair of points. 
    Data points are colored by the value of $k$. Other parameters are set to $A=0.1, R=40, F=1$. 
    Comparison to lines of slope $1$ (red) indicates that the run time increases at an approximately linear rate.	\label{fig:part-time1}}
\end{figure}

Next, we present an analysis of the run time and quality of the sampling on a DFN for varying sample sizes, variations of the parameter $k$, and different numbers of resampling attempts. 
All data is generated using the DFN shown in Fig.~\ref{fig:togetherzoomexp}. \jeffrey{is this correct?}\johannes{The histograms yes, since the run-time plots are for the 2d alg, they were only done on a few fractures from the figure, not the entire dfn.}
All samplings are performed on Fujitsu Laptop with 4 2.5 GHz intel Core i5 processors and 16GB of RAM \jeffrey{is this correct?}\johannes{fixed}
Different numbers of nodes were achieved by changing the parameter ${H}$, the minimum distance between points. 
All data is from independent samplings.
The plot in Figure \ref{fig:part-time1} shows the run time prior to the resampling process against the number of points sampled.
The color corresponds to the value of the parameter $k$, which controls the number of concurrent samples. 
We see an increase in run time with increasing values of $k$. 
The run times for samples with the same $k$ are positioned along straight lines of slope one, indicating a linear dependence of the total run time and the number of points sampled.
The red lines in the plot have a slope of one for reference. 

\begin{figure}
	\centering
	\input{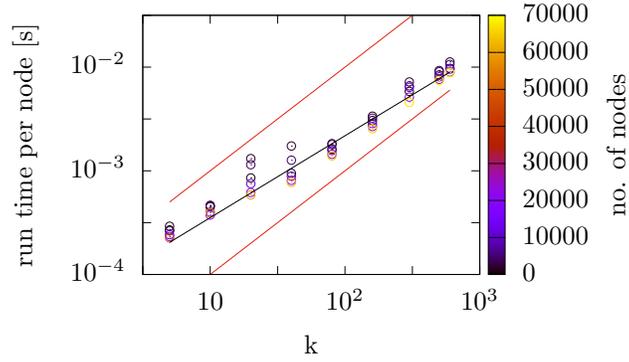}
	\caption{Run time per node of nMAPS plotted as a function of concurrently sampled points $k$ prior to the resampling process. 
    Data points generated using the same DFN.
    Different point densities generated by changing the minimum inhibition radius  $\frac{H}{2}$ between every pair of points. 
    Data points are colored by on the total number of points sampled. 
    Other parameters are set to $A=0.1, R=40, F=1$. 
    Linear fit (black) with slope $0.7(9)\pm 0.00(7)$. Comparison to lines of slope $1$(red) indicate sublinear behavior.}
	\label{fig:ktime}
\end{figure}

Figure \ref{fig:ktime} shows the relation between the parameter $k$ and the run time. Colors correspond to different numbers of points. 
As already established, the run time increases linearly with the number of points sampled.  The run time in terms of $k$ exhibits a slightly sublinear behavior. 
The linear fit (black) of the data on the log-log plot has a slope of $0.7(9)\pm 0.00(7)$. While this fitting error of $\approx 9\%$ is not insignificant, comparing the data to the two lines of slope $1$ (red) in the plot indicates that the run time does not increase more than linearly with $k$.

\begin{figure}
	\centering
	\input{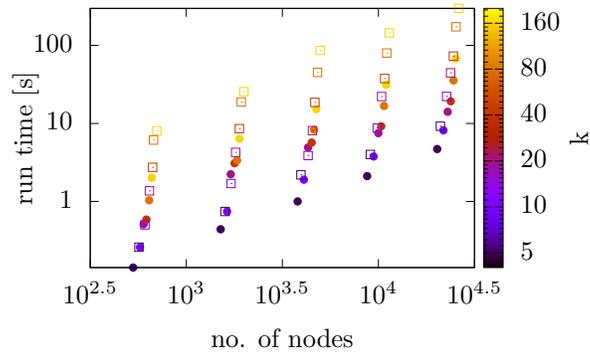}
	\caption{Comparison of run time for an implementation of \cite{Bridson2007FastPD,fastvar} (squares) and nMAPS algorithm (circles) plotted as a function of number of nodes and $k$ value. 
    Data points generated using the same DFN and different point densities using $H$. 
    Data points are colored depending on the value of $k$. 
    Other parameters are set to $A=0.1, R=40, F=1$. }
	\label{fig:nodestimecomp}
\end{figure}

Figure \ref{fig:nodestimecomp} shows a comparison of runtime for nMAPS variable-radii sampling with and without direct rejection of candidates using the background grid to find nearby points. 
Recall that the original method presented by \cite{Bridson2007FastPD} and \cite{fastvar} did not use direct rejection implemented in this manner. 
Thus, nMAPS without direct rejection of candidates is equivalent to the implementation of those algorithms.
Data points generated by nMAPS with direct rejection are represented by a filled circle, and data points generated without direct rejection are empty squares. 
All data points are colored by their $k$ value. 
nMAPS with direct rejection has a shorter run time for every pair of data points. 
The difference between the methods increases with larger values of $k$ because more candidates are rejected with larger values of $k$, but the direct rejection method does so without a distance calculation. 
For $k = 5$, the speed difference between the algorithms is slightly less than a factor of 2, whereas, for $k = 160$, the advantage has grown to about an order of magnitude. 
These plots highlight the performance gained by using the additional features found in nMAPS compared to previous methods.

A comparison of nMAPS with the method used in {\sc dfnWorks} \cite{hyman2015dfnWorks} to create a variable resolution point distribution is shown in Table \ref{tbl:runtime}. 
The method implemented in {\sc dfnWorks} uses an iterative Rivara refinement algorithm to generate the nodes of the mesh as described in Section~\ref{ssec:dfn}.
We consider three different DFN to characterize the difference between the methods. 
The first is the deterministic network of four ellipses shown in Fig.~\ref{fig:2d-together}.
The second is the network with fractures sampled from an exponential distribution containing 25 fractures shown in Fig.~\ref{fig:togetherzoomexp}.
The final network is composed of a single family of disc-shaped fractures whose fracture lengths are sampled from a truncated power-law with exponent 1.8, minimum length 1 m, maximum length 25 m within a cubic domain with sides of length 100 m. 
There are 8,417 fractures in this network.
In the first two examples, the algorithms were run on a MacBook Pro laptop with 8 2.9 GHz Intel Core i9 processors and 32GB of RAM.
The third example was run on a Linux server with 64 AMD Opteron(TM) Processor 6272 (1469.697 MHz) and 252GB of RAM.
Note that nMAPS is intrinsically parallelizable, and is implemented to mesh each fracture independently on a separate processor. 
The original meshing technique in {\sc dfnWorks} is parallelized in the same manner.  
The mesh resolution and setup were consistent between the methods so that the number of nodes in the final mesh is roughly the same.
In all cases, nMAPS was faster than the iterative method, and the speedup improved with the number of fractures. 
However, differences in network properties also likely play a role in the speedup, which is a feature that we do not explore in this study.

\begin{table}
    \begin{tabular}{|l|c|c|c|}
        \hline
        DFN Description & $\begin{matrix}
        	\text{Deterministic} \\ \text{Ellipses}
        \end{matrix}$ & $\begin{matrix}
        \text{Exponential} \\ \text{Distribution}
    \end{matrix}$ & $\begin{matrix}
    \text{Truncated} \\ \text{Power-law}
\end{matrix}$ \\ \hline
         Number of Fractures & 4 & 25 & 8417 \\
          Mesh Resolution & Uniform & Variable & Variable\\ 
          Number of Nodes & $\approx 1 \cdot 10^4$ & $\approx 1 \cdot 10^5$ & $\approx 3 \cdot 10^6$ \\ 
          Number of Processors & 4 & 8 & 32 \\
          Iterative - Run Time & 10.64 s  & 57.23 s & 47.47 m  \\
          nMAPS - Run Time &  4.46 s & 9.91 s &  6.47 m  \\
          Speed-up & 2.4 & 5.8 & 7.8 \\
    \hline
    \end{tabular}
\caption{Comparison of run time between previous iterative method with nMAPS for mesh generation. Number of nodes between methods differs slightly. In all cases, nMAPS outperforms the iterative method. \label{tbl:runtime}}
\end{table}


\subsection{Quality and resampling}
\begin{figure}
	\centering
	\input{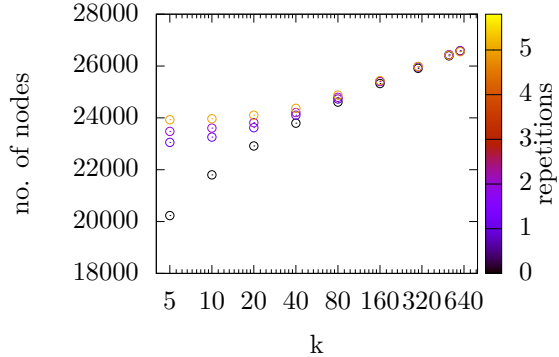}
	\caption{The total number of points accepted after resampling is plotted as a function of $k$ and colored by the number of resamplings. Other parameters are set to $A=0.1, R=40, F=1$.}
	\label{fig:k_nodes_rep}
\end{figure}

The maximality and density of nMAPS depends on the choice of $k$ and the number of times resampling is performed. 
Figure \ref{fig:k_nodes_rep} shows the total number of points accepted after a different number of resampling sweeps for various values of $k$.
Foremost, notice that the density of points increases with $k$. 
This rate of growth is largest for small values of $k$ and is lower at higher values. 
On the lower end of the $k$ scale, resampling increases the point density significantly, whereas there is barely a difference when $k>100$. 
The first resampling is particularly effective in adding additional nodes, whereas the difference between each additional resampling decreases thereafter. 
Given that resampling adds negligible run time to the original sampling process due to the efficient background mesh look-up in nMAPS, there is a trade-off between higher $k$ and more resampling sweeps. 
Adopting the latter method, i.e., a larger number of sweeps, can yield better performance in terms of obtaining a higher density with shorter run times. 
For example, a run at $k=5$ with a few resampling sweeps results in a density comparable to a run with more than $10$ times higher $k$ performed without resampling, but the former is significantly faster. 
Recall that run time increased linear with $k$.
Similar conclusions can be reached when looking at the quality of resulting triangulations rather than just the density of the Poisson-disk sampling.

\begin{figure}
	\centering
	\input{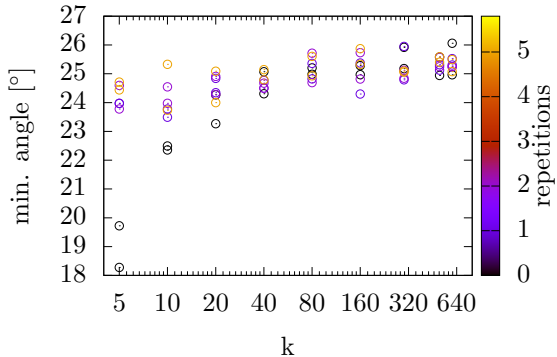}
	\caption{Smallest minimum angle of the triangulation produced using nMAPS. Points are colored by the number of resamplings. Increasing $k$ or the number of sweeps increases the minimum angle in the mesh. The latter of these requires much less time.}
	\label{fig:k_ang}
\end{figure}

Figure \ref{fig:k_ang} shows the smallest minimum angle in a triangulation using nMAPS for variable $k$ with different numbers of resampling sweeps. 
We can see for $k\gtrapprox 80$ this angle appears to be at around $25^\circ$ regardless of the number of repetitions. 
The theoretical bound for a maximal Poisson-disk sampling (with $r(\vx,\vy)=R(\vx,\vy )$) for the settings used to generate these data points is $27.04^\circ$.
Solving the the angle bounds from Lemma \ref{lemma:2} for $\varepsilon$ shows us that in this sampling $R(\vx,\vy)\lesssim (1+0.1)r(\vx,\vy)$.
Given that nMAPS is a stochastic method and identical inhibition/coverage radii are not exactly guaranteed, these results can be considered very good.
While the quality of triangulations for smaller $k$ values without resampling is significantly lower, it is worth emphasizing that a single repetition resolves this issue and produces triangulations with qualities on par with those produced using significantly higher $k$ values. 
Thus, nMAPS can be run using a single or low double-digit value of $k$, perform a single resampling sweep, and will produce a triangulation just as good as a higher value of $k$ would have produced but in a fraction of the time.

\subsection{Three-Dimensional Example}

\begin{figure*}
	\centering
\includegraphics[height=0.3\textwidth]{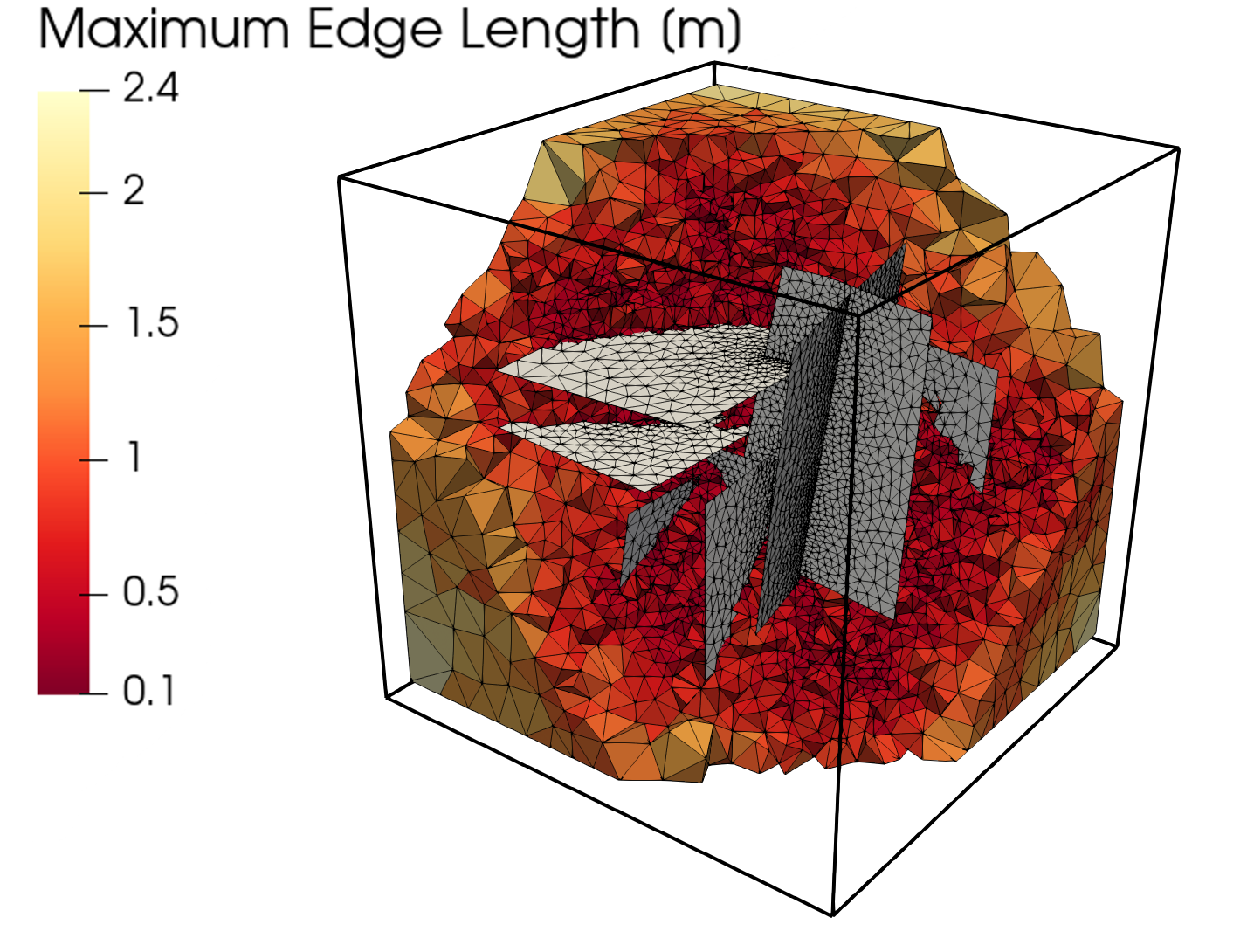}	\includegraphics[height=0.3\textwidth]{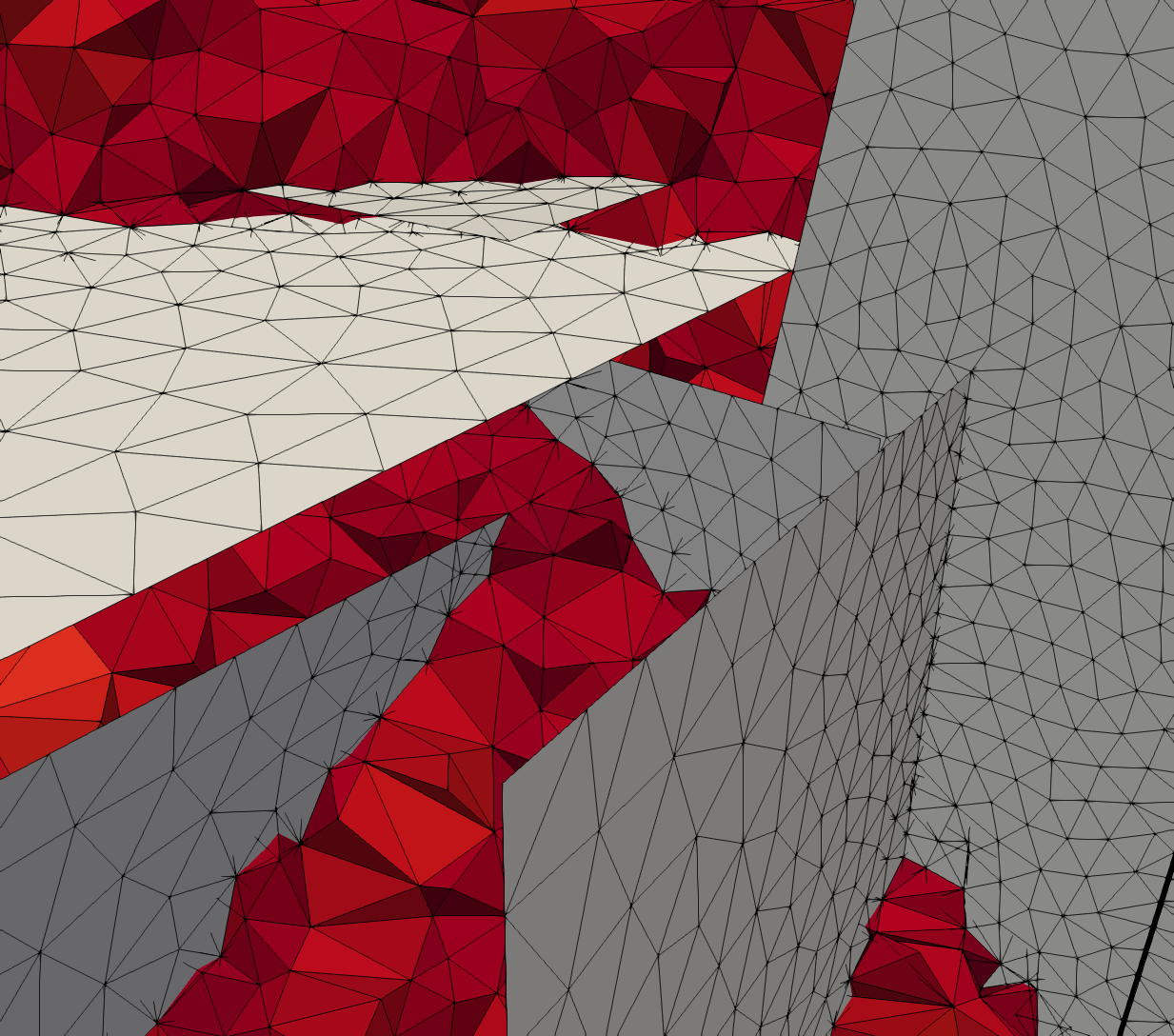}
	\caption{(Left) 3D Mesh produced by nMAPS of and around a seven fracture DFN. Parameters: $H=0.25$, $R=100$, $A=0.125$, and  $F=1$. (Right) Close up of the conforming mesh. Tetrahedra are colored according to their maximum edge length.}
	\label{fig:3d-mixzoomfin}
\end{figure*}

While nMAPS is primarily designed to optimize 2D sampling, we conclude with an example where these 2D samplings are combined with a 3D sampling of the surrounding volume to showcase how nMAPS can be used to produce high-quality 3D triangulations as well.
We use a simple network of seven square fractures for clarity in the visualization. 
The mesh produced by the 3D algorithm is shown in Fig.~\ref{fig:3d-mixzoomfin}. 
The tetrahedra are colored according to their maximum edge length to highlight how the point density depends on the distance from the DFN. 
Note that both the mesh of the DFN and the volume are variable resolution depending on the distance from the fracture intersections.

\begin{figure}
	\centering
	\begin{tikzpicture}
	\node (A) at (0,0)     {\input{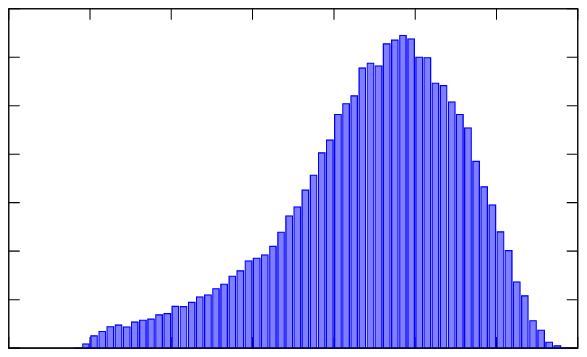}};
	\node (B) at (0,-5.5)  {\input{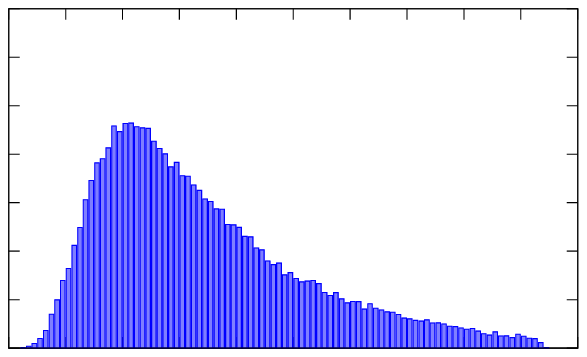}};
	\node (B) at (0,-11)   {\input{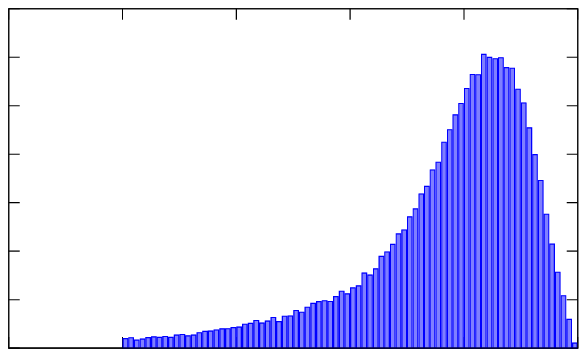}};
	\node (a) at (-4.,2.3) {(a)};
	\node (a) at (-4.,-3.2) {(b)};	
	\node (a) at (-4.,-8.7) {(c)};		
	\end{tikzpicture}
	
	\caption{Histograms of quality measures of the 3D mesh produced using nMAPS around a seven fracture DFN. Parameters: $H=0.01$, $R=40$, $A=0.1$, and $F=1$. (a): Minimum angle (all values $\ge 8^\circ$), (b): Maximum angle (all values $\le 165^\circ$), (c): Aspect Ratio (all values $\ge 0.2$)}
	\label{fig:3dhist}
\end{figure}

Histograms presented in Fig.~\ref{fig:3dhist} show the distribution of mesh quality measures of the tetrahedra in the 3D triangulation. 
Tetrahedra with either a dihedral angle of less than $8^\circ$ or an aspect ratio of less than $0.2$ are discarded before the sampling algorithm was restarted.  
Histogram (a) shows the distribution of the minimum dihedral angle of each tetrahedron. 
As expected, no dihedral angle below $8^\circ$ remains, while the vast majority exceeds values of $30^\circ$. 
Histogram (b) shows that despite not optimizing with respect to the maximum dihedral angle, none of these angles exceed $165^\circ$. 
Histogram (c) shows a sharp cut-off at $0.2$ in the distribution of aspect ratios, indicating that the aspect ratio is likely to have been the driving factor for a majority of the resamplings.  
This example was run through the sliver-removal resampling process $17$ times to obtain its triangulation quality. 
In each of these steps, only $200$ or fewer of $\approx 50000$ nodes were removed before the resampling. 
This low value of removed points indicates both the scarcity of slivers in samples generated through Algorithm 2 and that the vertices of these slivers can successfully be removed and replaced in a way that does not give rise to new slivers.

\section{Conclusions}\label{sec:conclusions}
We have presented the near-Maximal Algorithm for Poisson-disk Sampling (nMAPS) to produce a point distribution designed to generate high-quality variable resolution Delaunay triangulations in two and three dimensions.
We provided the theoretical basis on which nMAPS is built as well as details of its implementation. 
nMAPS uses efficient rejection and resampling techniques to achieve near maximality and linear run time scaling in the number of points accepted.
We demonstrated that nMAPS could successfully generate variable-radii Poisson-disk samples on polygonal regions, networks of polygons, and the surrounding volume they are embedded in.
Meshes generated using the point distributions produced by nMAPS show a quality nearly matching theoretical bounds for maximal Poisson-disk samplings, in which coverage and inhibition radii coincide.
It is worth noting that near maximality is reached for a coverage radius just slightly larger than the inhibition radius.
We determined that mesh quality produced by nMAPS is comparable to the method previously used in {\sc dfnWorks} but runs in a significantly shorter time.
Thus, nMAPS is significantly faster than the previous conforming variable mesh strategies due to our efficient rejection techniques that omit costly distance calculations.
It achieves mesh quality only marginally worse than what is theoretically possible. 
Moreover, it provides an iterative method where slivers in 3D volume meshes can be removed entirely from the domain within certain bounds.



It is worth mentioning that nMAPS is fast and simple to run in a parallel fashion in the context of mesh generation for DFNs, further improving the overall run time performance. 
nMAPS is intrinsically parallelizable by working on each fracture independently on a different processor, as was done in the performance comparison in Table~\ref{tbl:runtime}.
Based on the grid structure used to accept and reject candidates, 2D and 3D can also be further parallelized by dividing their domain into several pieces, which could be sampled individually on different processors while needing to communicate only cell information along the boundaries of the split domains. 
However, once these point distributions are produced, they all must reside on a single processor to connect them into a Delaunay mesh.

As a final comment, it's worth noting that nMAPS is not restricted to variable resolution DFN mesh generation for conforming numerical schemes.
As shown in Fig.~\ref{fig:2d-together}, it can be easily used to generate a uniform resolution mesh if desired. 
Moreover, nMAPS can be readily used to create meshes of fractures for non-conforming methods. 
One merely needs to skip the step in the algorithm where any point within the circumscribed circle of each segment of the discretized lines of intersection is removed from the point distribution. 
Omitting this step and performing a triangulation, not necessarily a Delaunay triangulation, will produce a mesh suitable for most non-conforming numerical schemes.
Finally, the general framework of nMAPS can be applied to efficiently produce arbitrary triangulations in two and three dimensions so long as the sampling can be constrained to within the domain. 

\section{Acknowledgments}
J.K. gratefully acknowledges support from the 2020 National Science Foundation
Mathematical Sciences Graduate Internship to conduct this research at Los Alamos National Laboratory.
J.D.H. and M.R.S. gratefully acknowledges support from the LANL LDRD program office Grant Numbers \#20180621ECR and \#20220019DR, the Department of Energy Basic Energy Sciences program (LANLE3W1), and the Spent Fuel and Waste Science and Technology Campaign, Office of Nuclear Energy, of the U.S. Department of Energy. 
M.R.S. would also like to thank support from the Center for Nonlinear Studies.
J.M.R. received support from DOE,  Contract No. DE-AC05-00OR22725. 
Los Alamos National Laboratory is operated by Triad National Security, LLC, for the National Nuclear Security Administration of U.S. Department of Energy (Contract No. 89233218CNA000001).
This work was prepared as an account of work sponsored by an agency of the United States Government.
The views and opinions of authors expressed herein do not necessarily state or reflect those of the United States Government or any agency thereof, its contractors or subcontractors.
LAUR \# LA-UR-21-24804.





\bibliography{references}
\bibliographystyle{plain}
\end{document}